\documentclass[a4paper]{amsart}
\usepackage{amsmath,amsthm,amsfonts}
\usepackage[all]{xy}
\newtheorem{thm}{Theorem}[subsection]
\newtheorem{thmx}{Theorem}[subsection]
\newtheorem{cor}[thm]{Corollary}
\newtheorem{lem}[thm]{Lemma}
\newtheorem{prop}[thm]{Proposition}
\theoremstyle{definition}
\newtheorem{defn}[thm]{Definition}
\theoremstyle{remark}
\newtheorem{rem}[thm]{Remark}
\newtheorem{ex}[thm]{Example}
\numberwithin{equation}{section}
\newcommand{\Rr}{\mathbb R}
\newcommand{\Zz}{\mathbb Z}

\newcommand{\set}[1]{\left\{#1\right\}}
\newcommand{\seq}[1]{\left\langle#1\right\rangle}
\newcommand{\eps}{\varepsilon}
\newcommand{\al}{\alpha}
\newcommand{\e}{\mathbf{e}}
\newcommand{\To}{\longrightarrow}

\newcommand{\D}{\mathcal{D}}
\newcommand{\F}{\mathcal{F}}
\newcommand{\Ho}{\mathcal{H}}

\newcommand{\X}{\mathcal{X}}
\newcommand{\Lie}{\mathcal{L}}

\newcommand{\Agerm}{\mathfrak{Aut}}
\newcommand{\Ogerm}{\mathfrak{Out}}
\newcommand{\Dgerm}{\mathfrak{Diff}}
\renewcommand{\gg}{\mathfrak{g}}

\newcommand{\gl}{\mathfrak{gl}}
\newcommand{\Exp}{\textrm{Exp}\,}
\newcommand{\Ker}{\textrm{Ker}~}
\newcommand{\Ad}{\textrm{Ad}\,}
\newcommand{\ad}{\textrm{ad}\,}
\newcommand{\Aut}{\textrm{Aut}\,}
\newcommand{\Out}{\textrm{Out}\,}
\newcommand{\Inn}{\textrm{Inn}\,}
\newcommand{\rank}{\textrm{rank}\,}
\newcommand{\tr}{\textrm{tr}\,}
\newcommand{\comment}[1]{}
\begin{document}

\title[Connections in Poisson Geometry]{Connections in Poisson Geometry I:\\
Holonomy and Invariants}
\author{Rui Loja Fernandes}
\address{Depart.~de Matem\'{a}tica, Instituto Superior T\'{e}cnico,
1049-001 Lisboa, PORTUGAL}
\email{rfern@math.ist.utl.pt}

\thanks{Supported in part by FCT grant PCEX/C/MAT/44/96 and PRAXIS XXI
through the Research Units Pluriannual Funding Program.}


\begin{abstract}
We discuss contravariant connections on Poisson manifolds. For
vector bundles, the corresponding operational notion of a
contravariant derivative had been introduced by I.~Vaisman. We show
that these connections play an important role in the study of
global properties of Poisson manifolds and we use them to define Poisson
holonomy and new invariants of Poisson manifolds.
\end{abstract}

\date{Jan.~31, 2000}
\maketitle


\section*{Introduction}

Let $M$ be a Poisson manifold and suppose that we require the existence
of a linear connection on $M$, compatible with the Poisson tensor
$\Pi$. Since parallel transport preserves the rank of the Poisson
tensor, the Poisson manifold must be regular in order for such
connection to exist. Therefore, the usual notion of a covariant
connection is not appropriate for the study of Poisson manifolds,
as some of the most interesting examples of Poisson manifolds are
non-regular. For non-regular Poisson manifolds the symplectic
foliation is singular and the dimension of the leaves varies, so one
can only hope to compare tangent spaces at different points of the
same symplectic leaf.

One possible way around this difficulty is to use families of
connections para\-me\-terized by the leaves. However, there are
examples showing that the  symplectic foliation can be wild, so the
space of leaves will not be easy to parameterize.

A much more efficient and direct approach, to be
introduced in this paper, is through the notion of a contravariant
connection, a concept that mimics for the case of Poisson manifolds
the usual notion of a covariant connection.

Assume we are given a principal bundle over a manifold $M$:
\[
\xymatrix{P \ar@(ur,dr)^G
\ar[d]_p\\ M}
\]
then a covariant connection $\Gamma$ on this principal bundle is
defined by a $G$-invariant horizontal distribution $u\mapsto H_u$
in $P$. Given a connection $\Gamma$, we have a notion of \emph{horizontal
lift}: $h(u,v)\in T_u P$ is the unique tangent vector to $H_u$
which projects to the vector $v\in T_{p(u)}M$. Conversely, the
horizontal lift $h$ defines the horizontal distribution
$H_u=\set{h(u,v): v\in T_{p(u)}M}$, so $h$ completely
determines the connection.

We shall define a
contravariant connection on a principal bundle over a Poisson
manifold by defining analogously the horizontal lift of cotangent
vectors. To formulate this notion, observe that $h$ is defined precisely for
pairs $(u,v)$ in $p^*TM$, the pullback bundle by $p$ of the tangent
bundle over $M$. Denote by $\widehat{p}:p^*TM\to TM$ the induced
bundle map so we have the commutative diagram
\[
\xymatrix{
p^*TM\ar[r]^{\widehat{p}}\ar[d]_{\widehat{\pi}}& TM \ar[d]^{\pi} \\
P\ar[r]_{p} &M }
\]
Then we can define a covariant connection to be a bundle map $h:p^*TM\to TP$,
such that:
\begin{enumerate}
\item[(CI)] $h$ is horizontal, i.~e., the following diagram commutes:
\[
\xymatrix{
p^*TM\ar[r]^{h}\ar[d]_{\widehat{p}}& TP \ar[d]^{p_*} \\
TM\ar[r]_{\text{id}} &TM }
\]
\item[(CII)] $h$ is $G$-invariant: $h(ua,v)=(R_a)_*h(u,v)$, for all $a\in G$;
\end{enumerate}

Assume now that $M$ is a Poisson manifold. According to a
general philosophical principle, in Poisson geometry sometimes the cotangent
bundle plays the role of the tangent bundle. Hence, we replace $TM$ by $T^*M$
in the diagrams above, whenever it makes sense. Thus we are lead to the
notion of
a \emph{contravariant connection} on a Poisson manifold: this is a bundle
map $h:p^*T^*M\to TP$, such that:
\begin{enumerate}
\item[(CI)$^*$] The following diagram commutes:
\[
\xymatrix{
p^*T^*M\ar[r]^{h}\ar[d]_{\widehat{p}}& TP \ar[d]^{p_*} \\ T^*
M\ar[r]_{\#} &TM }
\]
where $\#:T^*M\to TM$ is the bundle map induced by the Poisson tensor;
\item[(CII)$^*$] $h$ is $G$-invariant: $h(u a,\al)=(R_a)_*h(u,\al)$, for all
$a\in G$;
\end{enumerate}

Given a point $x$ in $M$ and a covector $\al\in T^*_x M$, the
vector $h(u,\al)\in T_u P$ will be called the \emph{horizontal
lift} of $\al$ to the point $u$ in the fiber over $x$. On any
fibration one can also consider \emph{generalized contravariant
connections} which satisfy only (CI)$^*$.

With such a definition at hand one can then develop the usual
concepts of pa\-rallelism, curvature, holonomy, geodesic, etc. In
particular, for a contravariant connection on a vector bundle
$p:E\to M$, one obtains in a way entirely analogous to the
covariant case, the notion of a
\emph{contravariant derivative} operator $D$: for each 1-form $\al$ on
$M$, $D_\al$ maps sections of $E$ to sections of $E$ and satisfies
\begin{enumerate}
\item[i)] $D_{\al+\beta}\phi=D_\al\phi+D_\beta\phi$;
\item[ii)] $D_\al(\phi+\psi)=D_\al\phi+D_\al\psi$;
\item[iii)] $D_{f\al}=fD_\al\phi$;
\item[iv)] $D_\al(f\phi)=fD_\al\phi+\#\al(f)\phi$;
\end{enumerate}
where $\al,\beta\in \Omega^1(M)$, $\phi$, $\psi$ are sections of
$E$, and $f\in C^{\infty}(M)$. Conversely, every such operator is
induced by a contravariant connection. Moreover, one can show that
there always exists a linear connection preserving the Poisson
tensor. In \cite{Vaisman:art:1} Vaisman introduces the notion of
contravariant derivative using i)-iv) as axioms.

In spite of its formal similarities with covariant connections,
there are striking differences in contravariant Poisson geometry.
For example, the holonomy of a connection may be non-discrete when
the connection is flat, contravariant connections cannot be pushed
back or forward, etc. However, just like in ordinary geometry,
contravariant connections are useful to study global properties of
Poisson manifolds.

Recall that the local structure of a Poisson manifold is given by
the Weinstein splitting theorem, also known as the generalized
Darboux theorem (see \cite{Weinstein:article:1}, Thm.~2.1). In a
neighborhood of a point, the Poisson structure splits as a direct
product of a symplectic structure and a Poisson structure which
vanishes at the point. So on the normal space to each symplectic
leaf we have a notion of
\emph{transverse Poisson structure}.

In global Poisson geometry one would like to understand the geometry and
topology of the symplectic foliation. Using generalized contravariant
connections we show that we have a notion of
\emph{Poisson holonomy} of the symplectic foliation, analogous to the
holonomy in the theory of regular foliations. The corresponding linear
holonomy coincides with the \emph{linear Poisson holonomy}
introduced by Ginzburg and Golubev in \cite{Ginzburg:article:1}.
The Poisson holonomy homomorphism is by Poisson automorphisms of
the transverse Poisson structure.

Poisson holonomy is not homotopy invariant, but factoring out the
inner Poisson automorphisms one obtains a notion of \emph{reduced
Poisson holonomy} invariant by homotopy, and we can prove the
following analogue of the Reeb stability theorem:

\begin{thmx}
Let $S$ be a compact, transversely stable leaf, with finite reduced
Poisson holonomy. Then $S$ is stable, i.~e., $S$ has arbitrarily
small neighborhoods which are invariant under all hamiltonian
automorphisms. Moreover, each symplectic leaf of $M$ near $S$ is a
bundle over $S$ whose fiber is a finite union of symplectic leaves
of the transverse Poisson structure.
\end{thmx}

We also discuss another related notion of holonomy, which we call
\emph{strict Poisson holonomy}, and which allows one to discuss global
splitting of an entire neighborhood of a symplectic leaf. The
corresponding stability theorem states that if $S$ has finite
strict Poisson holonomy then there is a neighborhood of $S$
which is Poisson covered by a product $\tilde{S}\times N$ where
$\tilde{S}$ is a finite cover of $S$.

Linear Poisson holonomy in turn can be discussed from the point of
view of linear contravariant connections and, for each symplectic
leaf, there is a notion of Bott contravariant connection. For a
non-regular Poisson manifold, we do not have a normal bundle (over
the whole of $M$) to the symplectic foliation. However, there is an
appropriate notion of a basic connection on $M$: these are linear
contravariant connections which preserve the Poisson tensor and
restrict in each leaf to the Bott contravariant connection.
Comparing a basic connection to a riemannian connection one is lead
to ``exotic" or secondary Poisson characteristic classes. These are
Poisson cohomology classes which give information on both the
Poisson geometry and the topology of the symplectic foliation of
$M$. In degree 1, this class actually coincides with the
\emph{modular class} of $M$. This invariant was discussed recently
by Weinstein in \cite{Weinstein:article:2}, where he shows that the
modular class is an obstruction to the existence of measures in $M$
invariant under the hamiltonian flows.

As a final note we remark that the most general setup for
contravariant connections is in the context of Lie algebroids.
Although we have omitted any references to Lie Algebroids, the
results discussed here should go through without any major changes, 
and this will be discussed elsewhere.

In a follow up to this paper (\cite{Fernandes:article:2}) we will
discuss invariant connections.
\vskip 10 pt

\begin{center}
\textsc{Acknowledgements}
\end{center}

This paper was certainly influenced by some remarks made by Alan
Weinstein after his talk at the Omega 99 Conference held in Lisbon,
which showed he had the complete picture on contravariant
connections on his mind. I also would like to thank my colleagues
Ana Cannas da Silva and Miguel Abreu for additional comments
and discussions.

\tableofcontents

\section{Contravariant Connections on Principal Bundles}

\subsection{Contravariant Cartan Calculus} On a Poisson
manifold there is a calculus on contravariant objects, analogous to
the usual Cartan calculus on differential forms. We recall here
some of the formulas and fix notation and conventions for later
use. Proofs of the results stated in this introductory paragraph
can be found in Vaisman's monograph \cite{Vaisman:book:1}.

Let $M$ be a Poisson manifold and denote by
$\Pi\in\X^2(M)$(\footnote{We denote by $\Omega^r(M)$ and $\X^r(M)$,
respectively, the spaces of differential $r$-forms and
$r$-multivector fields on a manifold $M$.}) the Poisson bivector
field, so the Poisson bracket on $M$ is given by
\begin{equation}
\set{f_1,f_2}=\Pi(df_1,df_2),\qquad f_1,f_2\in C^\infty(M).
\end{equation}
We also have a bundle map $\#:T^*M\To TM$ defined by
\begin{equation}
\beta(\#\al)=\Pi(\al,\beta),\qquad \al,\beta\in T^*M.
\end{equation}
On the space of differential 1-forms $\Omega^1(M)$ the Poisson
tensor induces a Lie bracket
\begin{equation}
\label{Lie:bracket:forms}
[\al,\beta]=\Lie_{\#\al}\beta-\Lie_{\#\beta}\al-d(\Pi(\al,\beta)),
\qquad \al,\beta\in \Omega^1(M),
\end{equation}
and for this Lie bracket and the usual Lie bracket on vector
fields, the map $\#:\Omega^1(M)\To\X^1(M)$ is a Lie algebra
homomorphism:
\begin{equation}
\label{eq:musical:homomorphism}
\#[\al,\beta]=[\#\al,\#\beta].
\end{equation}
We denote as usual by $X_f=\#(df)$ the hamiltonian vector field
associated with the function $f\in C^\infty(M)$, and we have
\begin{equation}
\label{eq:bracket:function}
[\al,f\beta]=f[\al,\beta]+\#\al(f)\beta=f[\al,\beta]-\left(i_{X_f}\al\right)
\beta.
\end{equation}

The existence of a Lie bracket on the space of 1-forms allows one
to mimic the algebraic definitions of $d$, $i_X$ and $\Lie_X$, to
obtain contravariant versions of these operators.

First, one defines the contravariant exterior differential
$\delta:\X^r(M)\To\X^{r+1}(M)$ by:
\begin{multline}
\delta Q(\al_0,\dots,\al_r)=\frac{1}{r+1}\sum_{k=0}^r (-1)^k\#\al_k(Q(\al_0,\dots,\widehat{\al}_k,\dots,\al_r)\\
+\frac{1}{r+1}\sum_{k<l}(-1)^{k+l}Q([\al_k,\al_l],\al_0,\dots,\widehat{\al}_k,\dots,\widehat{\al}_l,\dots,\al_r).
\end{multline}
where $\al_0,\dots,\al_r\in \Omega^1(M)$. This differential
satisfies:
\begin{align}
\delta^2(Q)&=0,\\
\label{eq:derivation}
\delta(Q_1\wedge Q_2)&=\delta Q_1 \wedge Q_2+(-1)^{\deg Q_1} Q_1\wedge \delta Q_2.
\end{align}
Moreover, if we extend the definition of $\#$ to forms of any
degree by setting
\begin{equation}
\#\lambda(\al_1,\dots,\al_r)=(-1)^r\lambda(\#\al_1,\dots,\#\al_r),
\end{equation}
we have
\begin{equation}
\delta(\#\lambda)=\#(d\lambda).
\end{equation}
The cohomology associated with $\delta$ is called the
\emph{Poisson cohomology} of $M$ and is denoted by $H^*_{\Pi}(M)$.
This relation shows that there is a homomorphism from de Rham cohomology
to Poisson cohomology $\#:H^*(M)\to H^*_{\Pi}(M)$, which in the
case of a symplectic manifold is an isomorphism.

Next, for each form $\al\in\Omega^1(M)$ there is an operator of contraction
by $\al$, denoted $i_\al:\X^r(M)\To\X^{r-1}(M)$, and an operator of Lie
derivative in the direction of $\al$, denoted
$\Lie_\al:\X^r(M)\To\X^r(M)$, given by
\begin{align}
(i_\al Q)(\al_1,\dots,\al_{r-1})&=Q(\al,\al_1,\dots,\al_{r-1}),\\
(\Lie_\al
Q)(\al_1,\dots,\al_r)&=\#\al(Q(\al_1,\dots,\al_r))-\sum_{k=1}^r
Q(\al_1,\dots,[\al,\al_k],\dots,\al_r).
\end{align}
We have formulas analogous to the usual formulas from Cartan
calculus:
\begin{align}
i_{[\al,\beta]}&=\Lie_\al i_\beta-i_\beta\Lie_\al,\\
\Lie_{[\al,\beta]}&=\Lie_\al \Lie_\beta-\Lie_\beta\Lie_\al,\\
\Lie_\al&=i_\al\delta+\delta i_\al,\\
\delta\Lie_\al&=\Lie_\al\delta.
\end{align}
In fact, the musical homomorphism relates these operators to the
usual ones, so for every 1-form $\al\in\Omega^1(M)$, every $r$-form
$\lambda\in\Omega^r(M)$ and every $r$-multivector field
$Q\in\X^r(M)$, one has:
\begin{align}
i_{\al}(\#\lambda)&=(-1)^r\#(i_{\#\al}\lambda),\\
\Lie_{\al}(\#\lambda)&=(-1)^r\#(\Lie_{\#\al}\lambda),\\
\Lie_{df} Q&=\Lie_{X_f}Q.
\end{align}

We can also extend $\Lie_\al$ to the exterior algebra $\Omega^*(M)$ by
setting
\begin{equation}
\Lie_\al \beta=[\al,\beta], \qquad \beta\in\Omega^1(M),
\end{equation}
and requiring $\Lie_\al$ to preserve type and act as a derivation.
Finally, we recall that the contravariant differential can also be
defined by
\begin{equation}
\delta Q=-[\Pi,Q]_s,
\end{equation}
where $[~,~]_s$ denotes the Schouten bracket.

\subsection{Contravariant Connections}
Let $P(M,G)$ be a smooth principal bundle over a Poisson manifold
$M$ with structure group $G$. We let $p:P\to M$ be the projection,
and for each $u\in P$ we denote by $G_u\subset T_u(P)$ the subspace
consisting of vectors tangent to the fiber through $u$. If we
denote by $p^*T^*M$ the pullback bundle, so there is a bundle map
$\widehat{p}:p^*T^*M \to T^*M$ which makes the following diagram
commutative
\[
\xymatrix{
p^*T^*M\ar[r]^{\widehat{p}}\ar[d]_{\widehat{\pi}}& T^*M \ar[d]^{\pi} \\
P\ar[r]_{p} &M }
\]
where on the vertical arrows we have the canonical projections.
Recalling that $p^*T^*M=\set{(u,\al)\in P\times
T^*M:p(u)=\pi(\al)}$, we see that we have a natural right
$G$-action on $p^*T^*M$ defined by $(u,\al)\cdot a\equiv(u a,\al)$,
if $a\in G$.

\begin{defn}
A \textsc{contravariant connection} $\Gamma$ in $P(M,G)$ is a
smooth bundle map $h:p^*T^*M\to TP$, such that:
\begin{enumerate}
\item[(CI)$^*$] The following diagram commutes:
\[
\xymatrix{
p^*T^*M\ar[r]^{h}\ar[d]_{\widehat{p}}& TP \ar[d]^{p_*} \\ T^*
M\ar[r]_{\#} &TM }
\]
\item[(CII)$^*$] $h$ is $G$-invariant: $h(u a,\al)=(R_a)_*h(u,\al)$, for all $a\in G$;
\end{enumerate}
\end{defn}

Given $(u,\al)\in p^*T^*M$, we call the vector $h(u,\al)\in T_u P$
the \emph{horizontal lift} of the 1-form $\al$ to $u$. The subspace
of $T_u P$ formed by all such horizontal vectors is denoted by
$\Ho_u$. The assignment $u\mapsto \Ho_u$ is a smooth, generalized,
distribution on $P$ called the \emph{horizontal distribution} of
the connection (by ``smooth'' we mean that for each point $u_0\in
P$ there exists a neighborhood $u_0\in U\subset P$ and smooth
vector fields $X_1,\dots,X_r$ in $U$, such that
$\Ho_u=\text{span}\set{X_1|_u,\dots,X_r|_u}$ for all $u\in U$).
Note that, as opposed to the covariant case, the rank of the
horizontal distribution will vary, and that this distribution does
not define the connection uniquely.

It follows from (CI)$^*$ in the definition of a contravariant
connection, that the horizontal spaces $\Ho_u$ project onto the
tangent space $T_x S$ to the symplectic leaf $S$ through $x=p(u)$.
In general, we have neither $T_u P=G_u+\Ho_u$ nor $G_u\cap
\Ho_u=\set{0}$. As usual, a vector $X\in T_u P$ will be called
\emph{vertical} (resp.~\emph{horizontal}), if it lies in $G_u$
(resp.~$\Ho_u$). If $M$ is not symplectic, a vector does not split
into a sum of an horizontal and a vertical component, so the usual
definitions of lift of curves, connection form, etc., do not make
sense in this context.

Later on, we shall need to consider \emph{generalized contravariant
connections}, by which mean that axiom (CII)$^*$ need not be
satisfied. Of course, such connections can be considered on any
fibration over a Poisson manifold.

\subsection{Connection Vector Fields}
If $\gg$ is the Lie algebra of $G$, we can express a contravariant
connection in $P$ by a family of $\gg$-valued vector fields, each
defined in an open subset of $M$. One should have in mind that, in
this theory, multivector fields play the role of differential
forms.

Henceforth, we use the following notation: We denote by $\set{U_j}$
an open cover of $M$, by $\psi_j:p^{-1}(U_j)\to U_j\times G$ a family
of trivializing isomorphisms, and by $\psi_{jk}:U_j\cap U_k \to G$ the
associated transition functions. For each $j$, we let $s_j:
U_j\to P$ be the section over $U_j$ defined by
$s_j(x)=\psi^{-1}_j(x,e)$, where $e\in G$ is the identity.

On each open set $U_j$ we define a $\gg$-valued vector field
$\Lambda_j$ as follows: if $\al\in\Omega^1(U_j)$, $x\in U_j$, and
$u=s_j(x)$, then
\[X_u=(s_j)_*\#\al_x-h(s_j(x),\al_x)\in T_u P\]
is a vertical vector since, by (CI)$^*$, we have:
\[ p_*X_u=p_*\cdot(s_j)_*\#\al_x-p_* h(s_j(x),\al_x)=\#\al_x-\#\al_x=0.\]
We let $\Lambda_j(\al)_x$ be the unique element $A\in\gg$ such that
$X_u=\sigma(A)_u$, which exists by (CII)$^*$. The $\set{\Lambda_j}$
are called the \emph{connection vector fields} of the contravariant
connection $\Gamma$.

In order to state the transformation rule for the connection vector
fields, it is convenient to introduce the following notation: if
$\phi:M\to N$ is a smooth map defined on a Poisson manifold $M$ its
\emph{contravariant differential} is the bundle map
$\delta \phi:T^*M\to TN$ defined by:
\begin{equation}
\delta \phi(\al_x)=d_x \phi\cdot\#\al_x, \qquad \al_x\in T^*_xM.
\end{equation}
If $N=\Rr$ this notation is consistent with the contravariant
differential introduced above, if we think of 0-vector fields as
functions.

\begin{prop}
The connection vector fields $\set{\Lambda_j}$ are related by
\begin{equation}
\label{eq:transform:connection}
\Lambda_k=\text{Ad}(\psi_{jk}^{-1})\Lambda_j+\psi_{jk}^{-1}\delta\psi_{jk},\quad \text{on }
U_j\cap U_k.
\end{equation}
Conversely, given a family of $\gg$-valued vector fields, each
defined in $U_j$, satisfying relations
(\ref{eq:transform:connection}), there is a unique contravariant
connection in $P(M,G)$ which gives rise to the $\set{\Lambda_j}$.
\end{prop}

\begin{proof}
Given a contravariant connection, define the vector fields
$\set{\Lambda_j}$ as above. If $U_j\cap U_k$ is non-empty, we have
$s_k(x)=s_j(x)\psi_{jk}(x)$, for all $x\in U_j\cap U_k$. If we set
$a=\psi_{jk}(x)\in G$, it follows from Leibniz rule that
\begin{equation}
\label{eq:section:translation}
{s_k}_*(X)=(R_{a})_*(s_j)_*(X)+\sigma((L_{a^{-1}})_*\cdot(\psi_{jk})_*X).
\end{equation}
If we compute both sides on $X=\#\al$, we obtain
\begin{align*}
\sigma(\Lambda_k(\al))_{ua}&={s_k}_*(\#\al)_{ua}-h(ua,\al)\\
                        &=(R_{a})_*(s_j)_*(\#\al)_u+\sigma((L_{a^{-1}})_*\cdot(\psi_{jk})_*\#\al)_u-(R_{a})_*h(u,\al)\\
                        &=(R_{a})_*\sigma(\Lambda_j(\al))_u+\sigma((L_{a^{-1}})_*\cdot(\psi_{jk})_*\#\al)_u\\
                        &=\sigma(\text{Ad}(\psi_{jk}^{-1})\Lambda_j(\al))_u+
                        \sigma(\psi_{jk}^{-1}\delta\psi_{jk}(\al))_u.
\end{align*}
as required.

Conversely, given a family of $\gg$-valued vector fields satisfying
relations (\ref{eq:transform:connection}), we define a
contravariant connection $\Gamma$ by letting the horizontal lift be
defined by
\begin{equation}
h(u,\al)={s_j}_*(\#\al)_u-\sigma(\Lambda_j(\al))_u,
\end{equation}
whenever $s_j$ is a section with $s_j(x)=u$. If $s_k$ is another
section with $s_k(x)=u$, it follows from
(\ref{eq:transform:connection}) and (\ref{eq:section:translation}),
with $\psi_{jk}(x)=a(x)=e$, that
\begin{align*}
{s_k}_*(\#\al)_u-\sigma(\Lambda_k(\al))_u&=
{s_j}_*(\#\al)_u+\sigma((\psi_{jk})_*\#\al)_u-\sigma(\Lambda_k(\al))_u\\
&={s_j}_*(\#\al)_u-\sigma(\Lambda_j(\al))_u,
\end{align*}
so this definition is independent of the section used. Conditions
(CI)$^*$ of the definition is easily verified. As for (CII)$^*$, we
note that if $\psi_{jk}(x)=a\in G$ is constant, then
$\Lambda_k=\text{Ad}(a^{-1})\Lambda_j$ and equation
(\ref{eq:section:translation}) gives
${s_k}_*(X)=(R_{a})_*(s_j)_*(X)$. Therefore, for any 1-form $\al$,
we find
\begin{align*}
h(ua,\al)&={s_k}_*(\#\al)_{ua}-\sigma(\Lambda_k(\al))_{ua}\\
&=(R_a)_*{s_j}_*(\#\al)_u-\sigma(Ad(a^{-1})\Lambda_j(\al))_{ua}=(R_a)_*h(u,\al),
\end{align*}
as wished.
\end{proof}

\subsection{Curvature}
For a contravariant connection $\Gamma$ with family of connection
vector fields $\set{\Lambda_j}$ we define a corresponding family of
\emph{curvature bivector fields} $\set{\Xi_j}$ by:
\begin{equation}
\label{eq:first:structure:equation}
\Xi_j=\delta \Lambda_j+\frac{1}{2}[\Lambda_j,\Lambda_j].
\end{equation}
Here, we are using the notation $[\xi,\zeta]$ for the $\gg$-valued
multivector field defined by
\[ [\xi,\zeta]=\sum_{a,b,c}C_{bc}^{a}\xi^b\wedge\zeta^c \e_a,\]
where $\xi=\sum_a \xi^{a} \e_a$ and $\zeta=\sum_a \xi^{a} \e_a$ are
$\gg$-valued multivector fields, relative to a basis $\set{\e_a}$
for $\gg$, and $C_{bc}^{a}$ are the structure constants of $\gg$
relative to the same basis.

\begin{prop}
The curvature bivector fields of a contravariant connection are
related by
\begin{equation}
\label{eq:transform:curvature}
\Xi_k=\text{Ad}(\psi_{jk}^{-1})\Xi_j,\quad \text{on }
U_j\cap U_k.
\end{equation}
Moreover, they satisfy the Bianchi identity:
\begin{equation}
\label{eq:Bianchi:1}
\delta \Xi_j+[\Lambda_j,\Xi_j]=0.
\end{equation}
\end{prop}

\begin{proof}
Set $\eta_{jk}=\psi_{jk}^{-1}\delta\psi_{jk}$. Then we have the
``Maurer-Cartan equations''
\begin{equation}
\delta \eta_{jk}=-\frac{1}{2}[\eta_{jk},\eta_{jk}].
\end{equation}
On the other hand, if $\Lambda_k$ and $\Lambda_j$ are related by
(\ref{eq:transform:connection}) we find, using
(\ref{eq:derivation}),
\begin{align*}
\delta\Lambda_k&=\delta(Ad(\psi^{-1}_{jk})\Lambda_j)+\delta\eta_{jk}\\
               &=+Ad(\psi^{-1}_{jk})\delta\Lambda_j
-\frac{1}{2}[\eta_{jk},Ad(\psi^{-1}_{jk})\Lambda_j]
+\frac{1}{2}[Ad(\psi^{-1}_{jk})\Lambda_j,\eta_{jk}]+\delta\eta_{jk}.
\end{align*}
Therefore, we have
\begin{align*}
\delta \Lambda_k+\frac{1}{2}[\Lambda_k,\Lambda_k]
&=Ad(\psi^{-1}_{jk})\delta\Lambda_j+
\frac{1}{2}[Ad(\psi^{-1}_{jk})\Lambda_j,Ad(\psi^{-1}_{jk})\Lambda_j]\\
&=Ad(\psi^{-1}_{jk})\left(\delta\Lambda_j+
\frac{1}{2}[\Lambda_j,\Lambda_j]\right).
\end{align*}
so (\ref{eq:transform:curvature}) holds.

Bianchi's identity (\ref{eq:Bianchi:1}) follows from
$\delta^2\Lambda_j=0$ and the derivation property
(\ref{eq:derivation}) of $\delta$.
\end{proof}

\begin{rem}
The structure equation (\ref{eq:first:structure:equation}) and the Bianchi
identity (\ref{eq:Bianchi:1}) show that
one should think of the operator $\delta+[\Lambda_j,\cdot]$ as a kind of
contravariant derivative acting on $\gg$-valued multivector fields.
This comment will be made precise later.
\end{rem}

It follows from (\ref{eq:transform:curvature}), that given 1-forms
$\al,\beta\in\Omega^1(M)$, we
can define a $\gg$-valued function $\Xi(\al,\beta)$ in $P$ by:
\[ \Xi(\al,\beta)_{s_j(x)}\equiv\Xi_j(\al,\beta).\]
$\Xi(\al,\beta)$ gives the following geometric interpretation of the
curvature:
Given a 1-form $\al\in\Omega^1(M)$, denote by
$h(\al)$ the horizontal lift of $\al$, so $h(\al)_u=h(u,\al)$ and
\[ u\mapsto \Ho_u=\set{h(\al)_u:\al\in\Omega^1(M)}\]
is the horizontal distribution.

\begin{prop}
\label{prop:curvature}
Let $\al,\beta\in\Omega^1(M)$. Then:
\begin{equation}
\label{eq:curvature:geometric}
[h(\al),h(\beta)]-h([\al,\beta])=-2\sigma(\Xi(\al,\beta)),
\end{equation}
\end{prop}

To prove the proposition we need the following lemma:

\begin{lem}
For any $\al,\beta\in\Omega^1(U_j)$
\begin{equation}
[h(\al),\sigma(\Lambda_j(\beta))]=-\sigma(\#\al(\Lambda_j(\beta))).
\end{equation}
\end{lem}

\begin{proof}
The flux of the vector field $\sigma(\Lambda_j(\beta))$ is
$\Phi_t(u)=u\exp(t\Lambda_j(\beta)(p(u)))$, so we have:
\[ [h(\al),\sigma(\Lambda_j(\beta))]_{u_0}=-\lim_{t\to 0}\frac{1}{t}
\left(
h(u_0,\al)-d\Phi_{-t}\cdot h(\Phi_t(u_0),\al)\right).\]
But:
\begin{align*}
d\Phi_{-t}\cdot h(\Phi_t(u_0),\al)&=
dR_{\exp(t\Lambda_j(\beta)(p(u)))}\cdot h(\Phi_t(u_0),\al)+
d\Psi\cdot h(\Phi_t(u_0),\al)\\
&=h(u_0,\al)+d\Psi\cdot h(\Phi_t(u_0),\al),
\end{align*}
where $\Psi(u)=u_0\exp(t\Lambda_j(\beta)(p(u)))$. Let
$s\mapsto\tilde{\gamma}(s,t)$ be the integral curve of $h(\al)$
through $\Phi_t(u_0)$. Then $s\mapsto
\gamma(s,t)=p(\tilde{\gamma}(s,t))$ is an integral curve of $\#\al$,
and we have:
\[ d\Psi\cdot h(\Phi_t(u_0),\al)=
\frac{d}{ds}\left.u_0\exp(t\Lambda_j(\beta)(\gamma(s,t)))\right|_{s=0}.\]
We conclude that
\begin{align*}
[h(\al),\sigma(\Lambda_j(\beta))]_{u_0}&=
-\frac{d}{dt}\left[\frac{d}{ds}\left.u_0\exp(t\Lambda_j(\beta)(\gamma(s,t)))\right|_{s=0}\right]_{t=0}\\
&=\frac{d}{ds}\left.\sigma(\Lambda_j(\beta)(\gamma(s,0)))_{u_0}\right|_{s=0}\\
&=\sigma(\#\al(\Lambda_j(\beta))_{p(u_0}))_{u_0},
\end{align*}
and the lemma follows.
\end{proof}

\begin{proof}[Proof of proposition \ref{prop:curvature}]
Over $U_j$ we have $(s_j)_*\#\al=\sigma(\Lambda_j(\al))+h(\al)$, so we
find:
\begin{align*}
[h(\al),h(\beta)]&=
(s_j)_*\#[\al,\beta]-[(s_j)_*\#\al,\sigma(\Lambda_j(\beta))]\\
&\qquad -[\sigma(\Lambda_j(\al)),(s_j)_*\#\beta]
+[\sigma(\Lambda_j(\al)),\sigma(\Lambda_j(\beta))]\\
&=h([\al,\beta])+\sigma(\Lambda_j([\al,\beta]))
-[h(\al),\sigma(\Lambda_j(\beta))]\\
&\qquad -[\sigma(\Lambda_j(\al)),h(\beta)]
-[\sigma(\Lambda_j(\al)),\sigma(\Lambda_j(\beta))])\\
&=h([\al,\beta])+\sigma(\Lambda_j([\al,\beta])-
\#\al(\Lambda_j(\beta))+\#\beta(\Lambda_j(\al))\\
&\qquad-\sigma([\Lambda_j(\al),\Lambda_j(\beta)])\\
&=h([\al,\beta])-\sigma(2\delta\Lambda_j(\al,\beta)+
[\Lambda_j,\Lambda_j](\al,\beta))\\
&=h([\al,\beta])-\sigma(2\Xi_j(\al,\beta)).
\end{align*}
\end{proof}

By a \emph{flat contravariant connection} we shall mean a
connection whose horizontal distribution is integrable.

\begin{prop}
A contravariant connection is flat iff its curvature bivector
fields vanish.
\end{prop}

\begin{proof}
By a result of Hermann \cite{Hermann:article:1}, a generalized
distribution associated with a vector subspace $\D\subset\X(M)$ is
integrable iff it is involutive and rank invariant. Taking
$\D=\set{h(df):f\in C^\infty(M)}$ so that
$\Ho_u=\set{X(u):X\in\D}$, proposition \ref{prop:curvature} shows
that $\D$ is involutive iff the curvature bivector fields vanish.
Hence, all it remains to show is that if the curvature vanishes and
$\gamma(t)$ is an integral curve of $h(df)$ then
$\dim\Ho_{\gamma(t)}$ is constant, for all small enough $t$.

Let $\tilde{\Phi}_t$ be the flow of $h(df)$ and let
$\Phi_t=p\circ\tilde{\Phi}_t$ be the flow of $\#df=X_f$.
If $\al\in\Omega^1(M)$ we claim that
\[ (\tilde{\Phi}_t)_* h(\al)=h(\Phi_{-t}^*\al), \]
for small enough $t$. In fact, the infinitesimal version of this
relation is
\[ [h(df),h(\al)]=h(\Lie_{X_f}\al)=h([df,\al]), \]
which by (\ref{eq:curvature:geometric}) holds, since we are assuming
that the curvature vanishes.

Therefore, the flow $\Phi_t$ gives an isomorphism between
$\Ho_{\gamma(0)}$ and $\Ho_{\gamma(t)}$, for small enough $t$, so $\D$
is rank invariant.
\end{proof}

\subsection{Parallelism and Holonomy}
Parallel displacement of fibers can be defined along curves
\emph{lying} on a symplectic leaf of $M$.

If $\gamma:[0,1]\to M$ is a smooth curve lying on a symplectic leaf
$S$, then $\gamma$ is also smooth as map $\gamma:[0,1]\to S$. This
follows from the existence of ``canonical coordinates'' for $M$ as
given by the generalized Darboux theorem. Also, by the same
theorem, we can choose a smooth family $t\mapsto\al(t)\in T^*M$ of
covectors such that $\#\al(t)=\dot{\gamma}(t)$. Following
\cite{Ginzburg:article:1}, we shall call the pair
$(\gamma(t),\al(t))$ a \emph{cotangent curve}.

\begin{prop}
Let $(\gamma(t),\al(t))$ be a cotangent curve. For any $u_0$
in $P$ with $p(u_0)=\gamma(0)$ there exists a unique horizontal
lift $\tilde{\gamma}:[0,1]\to P$, which satisfies the system
\begin{equation}
\label{eq:lift}
\left\{
\begin{array}{l}
\dot{\tilde{\gamma}}(t)=h(\tilde{\gamma}(t),\al(t)),\\
\\
\tilde{\gamma}(0)=u_0.
\end{array}
\right.
\end{equation}
\end{prop}

\begin{proof}
By standard results from the theory of o.d.e.'s with time dependent
coefficients, system (\ref{eq:lift}) has a unique maximal solution.
We claim that this solution exists for all $t\in [0,1]$.

By local triviality of the bundle we can find a curve $\bar{\gamma}:[0,1]\to P$
with $\bar{\gamma}(0)=u_0$ and $p(\bar{\gamma}(t))=\gamma(t)$. We look for a
curve $a(t)\in G$, such that $\tilde{\gamma}(t)=\bar{\gamma}(t)a(t)$
satisfies (\ref{eq:lift}). Differentiating, we have
\[\dot{\tilde{\gamma}}(t)=
\dot{\bar{\gamma}}(t)a(t)+\bar{\gamma}(t)\dot{a}(t).\]
We therefore require $a(t)$ to satisfy the equation
\[\dot{\bar{\gamma}}(t)a(t)+\bar{\gamma}(t)\dot{a}(t)=
h(\bar{\gamma}(t)a(t),\al(t)), \] or, equivalently,
\[ \bar{\gamma}(t)\dot{a}(t)a^{-1}(t)=
h(\bar{\gamma}(t),\al(t))-\dot{\bar{\gamma}}(t).\]
The right hand side of this equation belongs to $G_{\bar{\gamma}(t)}$ since
\[
p_*(h(\bar{\gamma}(t),\al(t))-\dot{\bar{\gamma}}(t))
=\#\al(t)-\frac{d}{dt}p(\bar{\gamma}(t))
=\#\al(t)-\dot{\gamma}(t)=0.
\]
Therefore, there exists some curve $A(t):[0,1]\to\gg$ such that
\[\bar{\gamma}(t)\dot{a}(t)a^{-1}(t)=\bar{\gamma}(t)A(t).\]
Since the initial value problem
\[\dot{a}(t)a^{-1}(t)=A(t), \qquad a(0)=e,\]
always has a solution, defined wherever $A(t)$ is defined, our claim
follows.
\end{proof}

Now using the proposition we can define parallel displacement of
the fibers along a cotangent curve $(\gamma(t),\al(t))$ in the
usual form: if $u_0\in p^{-1}(\gamma(0))$ we define
$\tau(u_0)=\tilde{\gamma}(1)$, where $\tilde{\gamma}(t)$ is the
unique horizontal lift of $(\gamma(t),\al(t))$ starting at $u_0$,
We obtain a map $\tau:p^{-1}(\gamma(0))\to p^{-1}(\gamma(1))$,
which will be called \emph{parallel displacement} of the fibers
along the cotangent curve $(\gamma(t),\al(t))$. It is clear, since
horizontal curves are mapped by $R_a$ to horizontal curves, that
parallel displacement commutes with the action of $G$:
\begin{equation}
\tau\circ R_a=R_a\circ\tau.
\end{equation}
Therefore, parallel displacement is an isomorphism between the
fibers.

If $x\in M$ lies in the symplectic leaf $S$, let $\Omega(S,x)$ be
the loop space of $S$ at $x$. Then for each cotangent loop
$(\gamma,\al)$, with $\gamma\in \Omega(S,x)$, parallel displacement
along $(\gamma,\al)$ gives a an isomorphism of the fiber
$p^{-1}(x)$ into itself. The set of all such isomorphisms forms the
holonomy group of $\Gamma$, with reference point $x$, and is
denoted $\Phi(x)$. Similarly, one has the restricted holonomy
group, with reference point $x$, denoted $\Phi^0(x)$, defined by
using cotangent loops in $S$ which are homotopic to the zero.

If $u\in p^{-1}(x)$ then we can also define the holonomy groups
$\Phi(u)$ and $\Phi^0(u)$. Just as in the covariant case, $\Phi(u)$
is the subgroup of $G$ consisting of those elements $a\in G$ such
that $u$ and $ua$ can be joined by an horizontal curve. We have
that $\Phi(u)$ is a Lie subgroup of $G$, whose connected component
of the identity is $\Phi^0(u)$, and we have isomorphisms
$\Phi(u)\simeq \Phi(x)$ and $\Phi(u)^0\simeq
\Phi(x)^0$.

If $x,y\in M$ belong to the same symplectic leave then the holonomy
groups $\Phi(x)$ and $\Phi(y)$ are isomorphic. This is because if
$u,v\in P$ are points such that, for some $a\in G$, there exists an
horizontal curve connecting $ua$ and $v$, then
$\Phi(v)=Ad(a^{-1})\Phi(u)$, so $\Phi(u)$ and $\Phi(v)$ are
conjugate in $G$. However, if $x,y\in M$ belong to different leaves
the holonomy groups $\Phi(x)$ and $\Phi(y)$ will be, in general,
non-isomorphic.

\begin{thm}
\label{thm:holonomy}
(Holonomy Theorem) Let $\Gamma$ be a contravariant
connection in $P(M,G)$, $u_0\in P$ and $s_j:U_j\to P$ a section
with $s_j(x_0)=p_0$. The Lie algebra of the holonomy group
$\Phi(u_0)\subset G$ is the ideal of $\gg$ spanned by all elements
of the form $\Xi_j(\al,\beta)_{x_0}+\Lambda_j(\gamma)_{x_0}$, where
$\al,\beta,\gamma\in T^*_{x_0}M$ are covectors with $\#\gamma=0$.
\end{thm}

\begin{proof}
It follows from the transformation rule
(\ref{eq:transform:connection}) for the connection vector fields,
that the subspace $\gg'\subset\gg$ spanned by all vectors of the
form $\Lambda_j(\gamma)_{x_0}$, with $\#\gamma=0$, is an ideal in
$\gg$. Similarly, it follows from the transformation rule
(\ref{eq:transform:curvature}) for the curvature bivector fields,
that the subspace $\gg''\subset\gg$ spanned by all vectors of the
form $\Xi(\al,\beta)_{x_0}$ is an ideal in $\gg$.

Let $P(u_0)$ be the set of points in $P$ that can be joined to
$u_0$ by a horizontal curve. We claim that the generalized
distribution $u\mapsto \Ho_u+\gg''_u$, where
$\gg''_u=\set{\sigma(A)_u:A\in\gg''}$, is integrable and that
$P(u_0)$ is the integral leaf through $u_0$. Assuming that this is
the case the proposition follows, for we have for any $A\in\gg$
\begin{align*}
A\in \text{Lie}(\Phi(u_0))\quad &\Longleftrightarrow \quad
a_t=\exp(tA)\in \Phi(u_0)\\ &\Longleftrightarrow \quad u_0 a_t\in
p^{-1}(p(u_0))\cap P(u_0)\\ &\Longleftrightarrow \quad\sigma(A)\in
G_{u_0}\cap T_{u_0}P(u_0)\\ &\Longleftrightarrow
\quad A\in\gg'+\gg''.
\end{align*}

The smooth distribution $u\mapsto \Ho_u+\gg''_u$ is integrable
becouse it is involutive and rank invariant. Let $G_{\Ho}$ be the
group of diffeomorphism generated by the horizontal vector fields
$h(\al)$. A theorem of Sussmann \cite{Sussmann:article:1}, shows
that the $G_{\Ho}$-invariant distribution $\D$ generated by $\Ho$
is integrable and that $P(u_0)$ is a leaf through $u_0$ of $\D$.
Therefore, the claim will follow if we can show that
$\D_u=\Ho_u+\gg''_u$. But, on one hand, $\D$ is involutive and
$\Ho\subset\D$, so we must have $\Ho_u+\gg''_u\subset\D_u$. On the
other hand, $\D$ is the smallest integrable distribution such that
$\Ho\subset\D$, so we must have $\D\subset \Ho_u+\gg''_u$.
\end{proof}

Note that the presence of the extra term $\gg'$ implies that a
connection can be flat and have non-discrete holonomy.

\subsection{Mappings of Connections}
Recall that a homomorphism $\phi:P(M,G)\to P'(M',G')$ of principal
bundles is a mapping of the total spaces $\phi:P\to P'$ such that
$\phi(ua)=\phi(u)\varphi(a)$, $u\in P$, $a\in G$, where $\varphi:G\to
G'$ is a Lie group homomorphism. We also have an induced map between the base
spaces, denoted here by the same letter: $\phi:M\to M'$. If this map
is a diffeomorphism and $s_j:U_j\to P$ is a local section of
$P(M,G)$ then $s'_j:\phi(U_j)\to P'$ defined by $s'_j=\phi\circ s_j\circ
\phi^{-1}$ is a local section of $P'(M',G')$.

\begin{prop}
\label{prop:connection:maps}
Let $M$ and $M'$ be Poisson manifolds and $\phi:P(M,G)\to P'(M',G')$ a
homomorphism such that the induced map $\phi:M\to M'$ is a Poisson
isomorphism. Given a contravariant connection $\Gamma$ in $P(M,G)$
there is a unique contravariant connection $\Gamma'$ in $P'(M',G')$
such that $\phi$ maps horizontal subspaces of $\Gamma$ to horizontal
subspaces of $\Gamma'$. The connection vector fields and the
curvature bivector fields of $\Gamma$ and $\Gamma'$ are related by:
\begin{equation}
\label{eq:trans:homomorphism:connection}
\Lambda'_j(\al)=\varphi_*\Lambda_j(\phi^*\al), \quad
\Xi'_j(\al,\beta)=\varphi_*\Xi_j(\phi^*\al, \phi^*\beta),\quad \al,
\beta\in\Omega^1(U'_j).
\end{equation}
If $u\in P$ and $u'=\phi(u)\in P'$, then $\varphi:G\to G'$ maps the holonomy
groups $\Phi(u)$ (resp.~$\Phi^0(u)$) onto $\Phi(u')$
(resp.~$\Phi^0(u')$).
\end{prop}

\begin{proof}
To define the connection $\Gamma'$, given $u'\in P'$ we choose
$u\in P$ and $a'\in G'$ such that $u'=\phi(u)a'$, and set
$h'(u',\al')=(R_{a'}\circ \phi)_* h(u,\phi^*\al')$. One checks that this
definition is independent of the choice of $u$ and $a'$.

If $b'\in G$, then $h'(u'b',\al')=(R_{a'b'}\circ \phi)_*
h(u,\phi^*\al')={R_{b'}}_*(R_{a'}\circ \phi)_*
h(u,\phi^*\al')={R_{b'}}_*h'(u',\al')$, hence $\Gamma'$ is invariant.
By invariance, we can now assume $\phi(u)=u'$, and we have:
\begin{align*}
p'_*h'(u',\al')&=p'_*\phi_*h(u,\phi^*\al')\\
                &=\phi_*p_*h(u,\phi^*\al')\\
                &=\phi_*\#\phi^*\al'=\#'\al',
\end{align*}
since $\phi:M\to M'$ is a Poisson map. Therefore, $\Gamma'$ is a
contravariant connection.

From the relation
\[ {s'_j}_*(\#'\al)=\phi_*{s_j}_*(\#\phi^*\al),\]
and the fact that the infinitesimal actions are related by
\[\sigma'(\varphi_*A)=\phi_*\sigma(A), \quad A\in\gg,\]
we obtain formulas (\ref{eq:trans:homomorphism:connection}) for the
connection vector fields. As for the curvature bivector fields we
have:
\begin{align*}
\Xi'_j(\al,\beta)&=\delta' \Lambda'_j(\al,\beta)+
               \frac{1}{2}[\Lambda'_j,\Lambda'_j](\al,\beta)\\
               &=\varphi_*\delta \Lambda_j(\phi^*\al, \phi^*\beta)+
               \frac{1}{2}[\varphi_*\Lambda_j,\varphi_*\Lambda_j](\phi^*\al, \phi^*\beta)\\
               &=\varphi_*\delta \Lambda_j(\phi^*\al, \phi^*\beta)+
               \frac{1}{2}\varphi_*[\Lambda_j,\Lambda_j](\phi^*\al, \phi^*\beta)
               =\varphi_*\Xi_j(\phi^*\al, \phi^*\beta),
\end{align*}
for any forms $\al,\beta\in\Omega^1(U'_j)$.

Finally, if $(\gamma',\al')$ is a cotangent loop at $x'=p'(u')$ lying in the
symplectic leaf through $x'$, then $(\gamma,\al)=(\phi^{-1}\circ
\gamma',\phi^*\al)$ is a cotangent loop at $x=p(u)$ lying in the
symplectic leaf through $x$. Therefore, if $\tilde{\gamma}$ is a
horizontal lift of $(\gamma,\al)$ then $\phi\circ \tilde{\gamma}$ is a
horizontal lift of $(\gamma',\al')$ and so the holonomy groups must be
related as stated.
\end{proof}

In the situation of the previous proposition we say that $\phi$ maps
the connection $\Gamma$ to the connection $\Gamma'$. There are two important
special cases to note:
\begin{enumerate}
\item[a)] if $P'(M',G')$ is a reduced sub-bundle of $P(M,G)$, so $M=M'$,
$\phi:M\to M$ is the identity map, and $h:G\to G'$ is a monomorphism,
we say the connection $\Gamma'$ is \emph{reducible} to the
connection $\Gamma$;
\item[b)] if $P'(M',G')=P(M,G)$, $M=M'$ and $\Gamma=\Gamma'$ we say
that the connection $\Gamma$ is \emph{invariant} by $\phi$, or simply
$\phi$-invariant. This means precisely that:
\begin{equation}
h(\phi(u),\al)=\phi_*h(u,\phi^*\al), \qquad \forall (\phi(u),\al)\in p^*T^*M;
\end{equation}
\end{enumerate}

For a general Poisson map it is not possible to pullback or
pushforward a contravariant connection, but there is still an
obvious definiton of \emph{mapping of connections}.

\subsection{Connections on Fiber Spaces}
\label{sec:fiberspaces}

If $G$ acts on the left on a manifold $F$ we shall denote by
$p_E:E(M,F,G,P)\to M$ the fiber bundle associated with $P(M,G)$
with standard fiber $F$.

Given a connection $\Gamma$ in $P(M,G)$ with associated horizontal
lift $h:p^*T^*M\to TP$, we define the induced horizontal lift
$h_E:p_E^*T^*M\to TE$ as follows: given $w\in E$ choose $(u,\xi)\in
P\times F$ which is mapped to $w$, and set
\begin{equation}
h_E(w,\al)=\xi_*h(u,\al),
\end{equation}
where we are identifying $\xi$ with the map $P\to E$ which sends
$u$ to the equivalence class of $(u,\xi)$. One can check easily
that this definition does not depend on the choice of $(u,\xi)$, so
we obtain a well defined bundle map $h_E:p_E^*T^*M\to TE$ which
makes the following diagram commute:
$$
\xymatrix{
p_E^*T^*M\ar[r]^{h_E}\ar[d]_{\widehat{p}_E}& TE \ar[d]^{p_{E*}} \\
T^*M\ar[r]_{\#} & TM
}
$$
As before, we can define horizontal and vertical vectors in $TE$,
horizontal lifts to $E$ of curves lying on symplectic leaves of
$M$, and parallel displacement of fibers of $E$. We shall call a
cross section $\sigma$ of $E$ over an open set $U\subset M$
\emph{parallel} if $\sigma_*(v)$ is horizontal for all tangent
vectors $v\in T_U M$.

\begin{thm} (Reduction Theorem)
\label{thm:reduction}
Let $P=P(M,G)$ be a principal fiber bundle over a Poisson manifold
$M$ with a contravariant connection $\Gamma$, and $H\subset G$ a
closed subgroup. There exists a one to one correspondence between
parallel cross sections $\sigma:M\to E(M,G/H,G,P)$ and sub-bundles
$Q(M,H)\subset P(M,G)$ such that $\Gamma$ is reducible to a
connection $\Gamma'$ in $Q$.
\end{thm}

\begin{proof}
Suppose we are given a parallel cross section $\sigma:M\to
E(M,G/H,G,P)$. Let $\pi:P\to E$ be the natural projection. Then we
define a sub-bundle $Q(M,H)$ by setting:
\[Q=\set{u\in P: \pi(u)=\sigma(p(u))}.\]
Given $u\in Q$ and $\al\in T^*_{p(u)}M$ let $(\gamma(t),\al(t))$ be
a cotangent curve with $\gamma(0)=p(u)$ and $\al(0)=\al$. The
horizontal lift $\tilde{\gamma}$ of this cotangent curve to $P$
satisfies $\mu(\tilde{\gamma}(t))=\sigma(\gamma(t))$, since
$\sigma$ is parallel. If follows that $h(u,\al)\in T_u Q$ for every
$u\in Q$, so $\Gamma$ is reducible to $Q$.

Conversely, suppose we are given a sub-bundle $Q(M,H)$ such that
$\Gamma$ is reducible to $Q$. Then we can define a section
$\sigma:M\to E(M,G/H,G,P)$ by setting $\sigma(x)=\pi(u)$, where
$u\in Q$ is any point satisfying $p(u)=x$. If $\tilde{\gamma}(t)$
is an horizontal curve in $P$ starting at $u\in Q$, then
$\tilde{\gamma}(t)\in Q$ since $\Gamma$ is reducible to $Q$. If
$\gamma(t)=p(\tilde{\gamma}(t))$, it follows that
$\mu(\tilde{\gamma}(t))$ is an horizontal lift of $\gamma$ to $E$
and that $\pi(\tilde{\gamma}(t))=\sigma(\gamma(t))$, so $\sigma$ is
flat.

\end{proof}

\subsection{Relationship to Ordinary Connections}
Let $M$ be a symplectic manifold and $\Gamma$ a contravariant
connection on $P(M,G)$ with horizontal lift $h:p^*T^*M\to TP$. Then
we have a bundle map $\tilde{h}:p^*TM\to TP$ defined by
\[ \tilde{h}(u,v)=h(u,\#^{-1}v),\qquad (u,v)\in p^*TM.\]
This map is obviously $G$-invariant and makes the following
diagram commute
\[
\xymatrix{
p^*TM\ar[r]^{\tilde{h}}\ar[d]_{\widehat{p}}& TP \ar[d]^{{p}_*} \\
TM\ar[r]_{\text{id}} &TM }
\]
It follows that $\tilde{h}$ is the horizontal lift of a covariant
connection on $M$. Let $\omega$ be the connection 1-form and let
$\Omega$ be the curvature 2-form of this connection. Also, given
trivialization isomorphisms $\set{\psi_j}$, inducing local sections
$\set{s_j}$, set $\omega_j=s_j^*\omega$ and $\Omega_j=s_j^*\Omega$.
Then it is clear from the definitions given above that the
connection vector fields $\set{\Lambda_j}$ and the curvature
bivector fields $\set{\Xi_j}$ are given by:
\begin{equation}
\Lambda_j=\#\omega_j,\qquad
\Xi_j=\#\Omega_j.
\end{equation}

For a general Poisson manifold with a contravariant connection
$\Gamma$ on $P(M,G)$ and horizontal lift $h:T^*M\to TP$, we say
that $\Gamma$ is \emph{induced by a covariant connection} if
\[ h(u,\al)=\tilde{h}(u,\#\al),\qquad (u,\al)\in p^*T^*M,\]
where $\tilde{h}:p^*TM\to TP$ is the horizontal lift of some
covariant connection on $M$. Note that in this case the lift $h$
satisfies:
\begin{equation}
\label{eq:F:connection}
\#\al=0\Longrightarrow h(u,\al)=0, \qquad (u,\al)\in p^*T^*M.
\end{equation}
This construction shows that there are always contravariant
connections on any principal bundle $P(M,G)$ over a Poisson manifold
$M$.

Not all connections satisfy property (\ref{eq:F:connection}), so we
set:

\begin{defn}
A contravariant connection $\Gamma$ on a principal bundle $P(M,G)$
is called a \textsc{$\F$-connection} if its horizontal lift satisfies
condition (\ref{eq:F:connection})
\end{defn}

Assume we have a contravariant $\F$-connection $\Gamma$ on
$P(M,G)$. If $i:S\hookrightarrow M$ is a symplectic leave, then on
the pull-back bundle $\tilde{p}:i^*P\to M$ we have an induced
connection $\Gamma_S$: on the total space $i^*P=\set{(y,u)\in
S\times P: i(y)=p(u)}$ we define the horizontal lift
$h_S:p_S^*T^*S\to T(i^*P)$ by setting
\begin{equation}
\label{eq:S:connection}
h_S((s,u),\al)=(p_*h(u,\beta),h(u,\beta)), \qquad (s,u)\in i^*P,
(u,\al)\in p^*T^*M,
\end{equation}
where $\beta\in T^*_{i(s)}M$ is such that $(d_s i)^*\beta=\al$, and
we are identifying $T(i^*P)=\set{(v,w)\in TS\times TP:v=p_*w}$. If
$(d_si)^*\beta'=(d_si)^*\beta$, then $\#\beta'=\#\beta$, so we get
the same result in (\ref{eq:S:connection}) and so $\Gamma$ is well
defined. $S$ being symplectic, the connection $\Gamma_S$ is induced
by a covariant connection on $i^*P$. Since the trivialization maps
$\psi_j:p^{-1}(U_j)\to U_j\times G$ induce trivialization maps
$\tilde{\psi}_j:\tilde{p}^{-1}(U_j\cap S)\to (U_j\cap S)\times G$
of the pull-back bundle $i^*P(M,G)$, writing
$\tilde{s}_j(y)=\tilde{\psi}^{-1}(y,e)$ for the associated sections,
we have:

\begin{prop}
\label{prop:F:connections}
Let $\Gamma$ be an $\F$-connection in $P(M,G)$. If $x\in M$
and $i:S\hookrightarrow M$ is the symplectic leaf
through $x$, denote by $\omega_S$ and $\Omega_S$ the connection
1-form and the curvature 2-form for the induced connection on
$i^*P(M,G)$. Also, let $\omega_j=\tilde{s}_j^*\omega_S$ and
$\Omega_j=\tilde{s}_j^*\Omega_S$. Then $\Lambda_j$ and $\Xi_j$ are
$i$-related to $\#\omega_j$ and $\#\Omega_j$:
\begin{equation}
\label{eq:relationship:connection}
i_*\#\omega_j=\Lambda_j,\qquad i_*\#\Omega_j=\Xi_j.
\end{equation}
\end{prop}

Therefore, a contravariant $\F$-connection in $P$ can be thought of
as a \emph{family} of ordinary connections over the symplectic
leaves of $M$. The (local) connection vector fields
$\set{\Lambda_j}$ and the (local) curvature bivector fields
$\set{\Xi_j}$ are obtained by gluing together the (local)
connection vector fields $\set{\#\omega_j}$ and the (local)
curvature bivector fields $\set{\#\Omega_j}$ of the connections on
the symplectic leaves of $M$.

For an $\F$-connection, horizontal lifts of cotangent curves
$(\gamma,\al)$ depend only on $\gamma$. Therefore, one has a well
determined notion of horizontal lift of a curve lying on a
symplectic leaf. It follows that for these connections, parallel
displacement can also be defined by first reducing to the pull-back
bundle over a symplectic leaf and then parallel displace the
fibers. Hence, the holonomy groups $\Phi(x)$ and $\Phi^0(x)$
coincide with the usual holonomy groups of the pull-back connection
on the symplectic leaf $S$ through $x$.

\subsection{Flat Connections}

Let $M$ be a Poisson manifold and $P(M,G)=M\times G$ the trivial
principal bundle. The \emph{canonical contravariant flat
connection} in $P(M,G)$ is defined by taking as horizontal lift
$h:p^*T^*M\to TP$ the map
\[h(u,\al)=(\#\al,0),\qquad (u,\al)\in p^*T^*M\]
where we identify $TP=TM\times TG$. This connection is a
$\F$-connection.

It is clear that a connection is the canonical flat connection iff
it is reducible to the unique contravariant connection in
$M\times{e}$, where $e\in G$ is the identity. For the canonical
flat connection and the natural trivialization the connection
vector field is $\Lambda=0$, and so the canonical flat connection
has zero curvature. Conversely, we have the following obvious
proposition:

\begin{prop}
For an $\F$-connection $\Gamma$ the following statements are
equivalent:
\begin{enumerate}
\item[i)] $\Gamma$ is flat;
\item[ii)] every point has neighborhood $U$ such
that the induced connection in $P|_U$ is isomorphic with the
canonical contravariant flat connection in $U\times G$;
\item[iii)] every point has neighborhood $U$ such that there
exists a parallel section $\sigma:U\to P$.
\end{enumerate}
Moreover, a flat $\F$-connection has discrete holonomy.
\end{prop}

If $\Gamma$ is not an $\F$-connection the conclusions of the
proposition, in general, do not hold.

\section{Linear Contravariant Connections}

\subsection{Contravariant Connections on a Vector Bundle}
Let $P(M,G)$ be a principal bundle over a Poisson manifold $M$ with
a contravariant connection $\Gamma$. Suppose that $G$ acts linearly
on a vector space $V$, so on the associated vector bundle
$E(M,V,G,P)$ we have the notion of parallel displacement of fibers
along cotangent curves $(\gamma,\al)$ (see section \ref{sec:fiberspaces}).

Given a section $\phi$ of $E$ defined along a cotangent curve
$(\gamma,\al)$, we define the \emph{contravariant derivative}
$D_{(\gamma,\al)}\phi$ to be the section
\begin{equation}
\label{eq:contra:derivative}
D_{(\gamma,\al)}\phi(t)=\lim_{h\to 0}\frac{1}{h}
\left[\tau^{t+h}_t(\phi(\gamma(t+h)))-\phi(\gamma(t))\right]
\end{equation}
where $\tau^{t+h}_t:p_E^{-1}(\gamma(t+h))\to p_E^{-1}(\gamma(t))$
denotes parallel transport of the fibers from $\gamma(t+h)$ to
$\gamma(t)$ along the cotangent curve $(\gamma,\al)$.

\begin{prop}
\label{prop:derivative:properties}
Let $\phi$ and $\psi$ be sections of $E$ and $f$ a function on $M$
defined along $\gamma$. Then
\begin{enumerate}
\item[i)] $D_{(\gamma,\al)}(\phi+\psi)=D_{(\gamma,\al)}\phi+D_{(\gamma,\al)}\psi$;
\item[ii)] $D_{(\gamma,\al)}(f\phi)=(f\circ\gamma)D_{(\gamma,\al)}\phi+\dot{\gamma}(f)(\phi\circ\gamma)$;
\end{enumerate}
\end{prop}

\begin{proof}
i) is obvious from the definition. On the other hand, we have
\[\tau^{t+h}_t(f(\gamma(t+h))\phi(\gamma(t+h)))=
f(\gamma(t+h))\tau^{t+h}_t(\phi(\gamma(t+h))),\] and ii) follows by
the Leibniz rule.
\end{proof}

Now let $\al\in T_x^*M$ be a covector and $\phi$ a cross section of
$E$ defined in a neighborhood of $x$. The contravariant derivative
$D_\al\phi$ of $\phi$ in the direction of $\al$ is defined as
follows: choose a cotangent curve $(\gamma(t),\al(t))$ defined for
$t\in(-\eps,\eps)$, and such that $\gamma(0)=x$ and $\al(0)=\al$.
Then we set:
\begin{equation}
\label{eq:contra:derivative:0}
D_\al\phi=D_{(\gamma,\al)}\phi(0).
\end{equation}
It is easy to see that $D_\al\phi$ is independent of the choice of
cotangent curve. Clearly, a cross section $\phi$ of $E$ defined on
an open set $U\subset M$ is flat iff $D_\al\phi=0$ for all $\al\in
T_x M$, $x\in M$.

\begin{prop}\label{prop:derivative:properties:1}
Let $\al,\beta\in T_x^*M$, $\phi$ and $\psi$ cross sections of
$E$ defined in a neighborhood $U$ of $x$. Then
\begin{enumerate}
\item[i)] $D_{\al+\beta}\phi=D_\al\phi+D_\beta\phi$;
\item[ii)] $D_\al(\phi+\psi)=D_\al\phi+D_\al\psi$;
\item[iii)] $D_{c\al}=cD_\al\phi$, for any scalar $c$;
\item[iv)] $D_\al(f\phi)=f(x)D_\al\phi+\#\al(f)\phi(x)$, for any function $f\in C^{\infty}(U)$;
\end{enumerate}
\end{prop}

\begin{proof}
iii) is obvious, while ii) and iv) follow from proposition
\ref{prop:derivative:properties}. To prove i) observe that any
section $\phi$ of $E$, defined in a open set $U$, can be identified
with a function $F:p^{-1}(U)\to V$ by letting
\[ F(u)=u^{-1}(\phi(p(u))), \qquad u\in p^{-1}(U),\]
where we view $u\in P$ as a linear isomorphism $u:V\to p^{-1}_E(u)$.
Then, as in the covariant case, we find
\[ D_\al\phi=u(h(u,\al)\cdot F).\]
From this expression for the contravariant derivative, i) follows
immediately.
\end{proof}

Now let $\al\in\Omega^1(M)$ be a 1-form and $\phi$ a section of
$E$. We define the contravariant derivative $D_\al\phi$ to be the
section of $E$ given by:
\begin{equation}
\label{eq:contra:derivative:00}
D_\al\phi(x)=D_{\al_x}\phi.
\end{equation}

\begin{prop}
\label{prop:derivative:properties:2}
Let $\al,\beta\in \Omega^1(M)$, $\phi$ and $\psi$ cross sections of
$E$, and $f\in C^{\infty}(M)$. Then
\begin{enumerate}
\item[i)] $D_{\al+\beta}\phi=D_\al\phi+D_\beta\phi$;
\item[ii)] $D_\al(\phi+\psi)=D_\al\phi+D_\al\psi$;
\item[iii)] $D_{f\al}=fD_\al\phi$;
\item[iv)] $D_\al(f\phi)=fD_\al\phi+\#\al(f)\phi$;
\end{enumerate}
\end{prop}

\begin{proof}
From proposition \ref{prop:derivative:properties:1} we obtain
immediately that i)-iv) hold.
\end{proof}

It is also true that the contravariant derivative uniquely
determines the connection. The proof of the following proposition
is similar to the covariant case and so it will be omitted.

\begin{prop}
Suppose for each 1-form $\al\in\Omega^1(M)$ there is a linear operator
$D_\al$ acting on sections of $E$ and satisfying i)-iv) of proposition
\ref{prop:derivative:properties:2}.
Then there exists a unique contravariant connection $\Gamma$ on the
associated principal bundle $P(M,G)$ whose induced contravariant
derivative on $E$ is $D$.
\end{prop}

In the case where the contravariant connection is induced by a
covariant connection, the contravariant derivative $D$ and the
covariant derivative $\nabla$ are related by
\begin{equation}
D_\al =\nabla_{\#\al}.
\end{equation}

On the other hand, $\F$-connections can be characterized by the
condition:
\begin{equation}
\#\al=0\ \Longrightarrow\ D_\al=0, \qquad \forall \al\in T^*(M).
\end{equation}
Moreover, by proposition \ref{prop:F:connections}, for an
$\F$-connection, on each symplectic leaf $i:S\hookrightarrow M$
there is a covariant connection on the pullback bundle $i^*P$,
inducing a covariant derivative $\nabla$ on $i^*E$, with the
following property: if $\psi$ is any cross section of $E$, then
\begin{equation}
i^*D_\al \psi =\nabla_{\#i^*\al}i^*\psi,
\end{equation}
where $i^*\psi$ denotes the section of the pullback bundle $i^*E$
induced by $\psi$.

\comment{Finally, if the structure group $G=GL(n)$ we can choose the
standard basis for $\gg=\mathfrak{gl}(m)$ formed by the elementary
matrices $E_i^j$, with 1 in the $j^{\textbf{th}}$-row and
$i^{\textbf{th}}$-column and 0 elsewhere.  Then the connection
vector fields and the curvature bivector fields can be written in
the form:
\begin{equation}
\Lambda_k=\sum_{i,j} {\Lambda}_{kj}^i E^j_i,\qquad
\Xi_k=\sum_{i,j} \Xi_{kj}^i E^j_i.
\end{equation}
In this notation, these multivector fields are related by:
\begin{equation}
\delta {\Lambda}_{kj}^i=\Xi_{kj}^i-
\sum_{l=1}^m {\Lambda}_{kl}^i\wedge {\Lambda}_{kj}^l.
\end{equation}
}

\subsection{Linear Contravariant Connections}

A \emph{linear contravariant connection} is a contravariant connection on the
coframe bundle $P=F^*(M)$ over $M$, so $G=GL(m)$ where $m=\dim M$.
If $u=(\al_1,\dots,\al_m)\in F^*(M)$ is a coframe, we can view $u$ as a
linear isomorphism
$u:{(\Rr^m)}^*\to T^*_{p(u)}M$ by setting
\[ u(\xi)(v)=\xi(\al_1(v),\dots,\al_m(v)), \qquad v\in T_{p(u)}M,\
\xi\in{(\Rr^m)}^*.\]
We define the \emph{canonical vector fields} $\theta_j$ on an open set $U_j$,
with trivializing isomorphism $\psi_j:p^{-1}(U_j)\to U_j\times G$, and
associated section $s_j(x)=\psi_j^{-1}(x,e)$, to be the ${(\Rr^m)}^*$-valued
vector fields defined by
\begin{equation}
\label{eq:canonical:fields}
\theta_j(\al)_x=s_j(x)^{-1}(\al), \qquad x\in U_j.
\end{equation}
These allows us to define the \emph{torsion bivector fields}
$\Theta_j$ to be the ${(\Rr^m)}^*$-valued bivector fields given by
\begin{equation}
\label{eq:torsion:fields}
\Theta_j(\al,\beta)=\delta\theta_j(\al,\beta)+
\Lambda_j(\al)\cdot\theta(\beta)-\Lambda_j(\beta)\cdot\theta_j(\al).
\end{equation}

\begin{prop}
The canonical vector fields and the torsion bivector fields of a linear
contravariant connection are related by
\begin{eqnarray}
\label{eq:transform:canonical}
  \theta_k&=\psi_{jk}^{-1}\cdot\theta_j,\\
\label{eq:transform:torsion}
  \Theta_k&=\psi_{jk}^{-1}\cdot\Theta_j.
\end{eqnarray}
Moreover, they satisfy the Bianchi identity
\begin{equation}
\label{eq:Bianchi:2}
\delta \Theta_j(\al,\beta,\gamma)=
\bigodot_{\al,\beta,\gamma}\delta\Lambda_j(\al,\beta)\cdot\theta_j(\gamma)-
\bigodot_{\al,\beta,\gamma}\Lambda_j(\al)\cdot\delta\theta_j(\beta,\gamma).
\end{equation}
where the symbol $\bigodot$ denotes cyclic sum over the subscripts.
\end{prop}

\begin{proof}
Relation (\ref{eq:transform:canonical}) follows immediately from the
definition of the canonical vector fields. To prove
(\ref{eq:transform:torsion}), we take the contravariant differential of
(\ref{eq:transform:canonical}):
\[\delta\theta_k(\al,\beta)=\psi_{jk}^{-1}\cdot\delta\theta_j(\al,\beta)
-\psi_{jk}^{-1}\delta\psi_{jk}(\al)\psi_{jk}^{-1}\cdot\theta_j(\beta)+
\psi_{jk}^{-1}\delta\psi_{jk}(\beta)\psi_{jk}^{-1}\cdot\theta_j(\al).\]
From the transformation rule (\ref{eq:transform:connection}) for
the connection vector fields, we find
\[ \Lambda_k(\al)\cdot\theta_k(\beta)=
\psi_{jk}^{-1}\Lambda_j(\al)\cdot\theta_j(\beta)+
\psi_{jk}^{-1}\delta\psi_{jk}(\al)\psi_{jk}^{-1}\cdot\theta_j(\beta).\]
Therefore, we compute:
\begin{align*}
\Theta_k(\al,\beta)&=\delta\theta_k(\al,\beta)+
\Lambda_k(\al)\cdot\theta(\beta)-\Lambda_k(\beta)\theta_j(\al)\\
&=\psi_{jk}^{-1}\delta\theta_j(\al,\beta)+
\psi_{jk}^{-1}\Lambda_j(\al)\cdot\theta_j(\beta)-
\psi_{jk}^{-1}\Lambda_j(\beta)\cdot\theta_j(\al)=
\psi_{jk}^{-1}\cdot\Theta_j(\al,\beta).
\end{align*}

The Bianchi identity follows from taking the contravariant
differential of (\ref{eq:torsion:fields}).
\end{proof}

For the standard contragradient action of $G=GL(m)$ on $F=(\Rr^m)^*$, the
bundle associated with the coframe bundle $P=F^*(M)$ is the cotangent bundle
$T^*M=E(M,F,G,P)$. Sections of $T^*(M)$ are just differential 1-forms and so
the contravariant derivative associates to each 1-form $\al$ a linear operator
$D_\al:\Omega^1(M)\to\Omega^1(M)$ such that:
\begin{align}
\label{eq:contra:derivative:1}
&D_{f_1\al_1+f_2\al_2}=f_1 D_{\al_1}+f_1 D_{\al_1}, \qquad
\text{for all } f_i\in C^\infty(M),\ \al_i\in\Omega^1(M),\\
\label{eq:contra:derivative:2}
&D_\al (f\beta)=f D_\al \beta+\#\al(f)\beta,\qquad \text{for all }
f\in C^\infty(M),\
\al,\beta\in\Omega^1(M).
\end{align}

One can also consider other associated vector bundles to $F^*(M)$ which
lead, just us in the covariant case, to contravariant derivatives of any
tensor fields over $M$. For example, if $X$ is a vector field, then
$D_\al X$ is the contravariant derivative of $X$ along the 1-form $\al$.
It is completely characterized by the relation
\begin{equation}
\seq{D_\al X,\beta}=\#\al(\seq{X,\beta})-\seq{X,D_\al \beta},
\end{equation}
which holds for every 1-form $\beta\in\Omega^1(M)$. One has similar formulas
for the contravariant derivative of any tensor field on $M$.

Local coordinate expressions for linear contravariant connections can
be obtained in a way similar to the covariant case. Let
$(x^1,\dots,x^m)$ be local coordinates on a neighborhood $U$ in
$M$. Then we define Christoffel symbols $\Gamma^{ij}_k$ by
\begin{equation}
D_{dx^i}dx^j =\Gamma^{ij}_k dx^k.
\end{equation}
It is easy to see that under a change of coordinates these symbols
transform according to
\begin{equation}
\label{eq:symbols:transform}
{\tilde{\Gamma}}^{lm}_n=
\frac{\partial y^l}{\partial x^i}\frac{\partial y^m}{\partial x^j}
\frac{\partial x^k}{\partial y^n}\Gamma^{ij}_k
+
\frac{\partial y^l}{\partial x^i}\frac{\partial^2 y^m}{\partial x^j \partial x^k}
\frac{\partial x^j}{\partial y^n}\pi^{ik},
\end{equation}
where $\pi^{ik}$ are the components of the Poisson tensor.
Conversely, given a family of symbols that transform according to
this rule under a change of coordinates, we obtain a well defined
contravariant derivative/connection on $M$.

Using these symbols, it is easy to get the local coordinates
expressions for the contravariant derivatives: given a 1-form
$\al=\al_i dx^i$ and a tensor field $K$, of type $(r,s)$, with
components $K^{i_1\dots i_r}_{j_1\dots j_s}$, we have
\begin{multline}
\label{eq:derivative}
(D_\al K)^{i_1\dots i_r}_{j_1\dots j_s}=
\pi^{kl}\al_k\frac{\partial K^{i_1\dots i_r}_{j_1\dots j_s}}{\partial x^l}
-\sum_{a=1}^r\left(\Gamma^{k i_a}_l\al_k K^{i_1\dots l\dots
i_r}_{j_1\dots j_s}\right)\\ +\sum_{b=1}^s\left(\Gamma^{k
l}_{j_b}\al_k K^{i_1\dots i_r}_{j_1\dots l\dots j_s}\right).
\end{multline}

Given a tensor field $K$ of type $(r,s)$ we shall write, as in the
covariant case, $DK$ for the tensor field of type $(r+1,s)$ such
that
\begin{equation}
(DK)^{i_1\dots i_r k}_{j_1\dots j_s}=(D_{dx^k} K)^{i_1\dots
i_r}_{j_1\dots j_s}.
\end{equation}
A tensor field $K$ on $M$ is \emph{parallel} iff $DK=0$.

\subsection{Curvature and Torsion Tensor Fields}

For a linear contravariant connection on a Poisson manifold $M$ we define
the \emph{torsion tensor field} $T$ and the \emph{curvature tensor
field} $R$, respectively, to be the tensor fields of types $(2,1)$
and $(3,1)$ given by
\begin{align}
\label{torsion}
T(\al,\beta)&=s_j(x)(\Theta_j(\al,\beta),\\
\label{curvature}
R(\al,\beta)\gamma&=s_j(x)\left[\Xi^*_j(\al,\beta)\cdot
s_j^{-1}(x)(\gamma)\right].
\end{align}
where $x\in U_j$, $\al,\beta,\gamma\in T_x^*(M)$, and we are denoting
by $\Xi^*_j(\al,\beta)$ the endomorphism of $\mathfrak{gl}(m)$ dual
to $\Xi_j(\al,\beta)$. Note that if $x\in U_j\cap U_k$ and
$s_k(x)=\psi_{jk}(x)s_k(x)$ we obtain the same values in formulas
(\ref{torsion}) and (\ref{curvature}), so these really define
tensor fields on all of $M$. These tensor fields can be easily
expressed in terms of contravariant derivatives:
\begin{prop}
In terms of contravariant differentiation, the torsion $T$ and the
curvature $R$ can be expressed as follows:
\begin{align}
\label{eq:torsion:derivative}
T(\al,\beta)&=D_\al \beta-D_\beta \al - [\al,\beta],\\
\label{eq:curvature:derivative}
R(\al,\beta)\gamma&=D_\al D_\beta \gamma-D_\beta D_\al \gamma
-D_{[\al,\beta]}\gamma.
\end{align}
Moreover, the Bianchi identities (\ref{eq:Bianchi:2}) and
(\ref{eq:Bianchi:1}) can also be expressed as
\begin{align}
&\bigodot_{\al,\beta,\gamma}\left(D_\al R(\beta,\gamma)+
R(T(\al,\beta),\gamma)\right)=0,\\
&\bigodot_{\al,\beta,\gamma}\left(R(\al,\beta)\gamma-
T(T(\al,\beta),\gamma)-D_\al T(\beta,\gamma)\right)=0.
\end{align}
\end{prop}

From formulas (\ref{eq:torsion:derivative}) and
(\ref{eq:curvature:derivative}), we obtain immediately the
following local coordinates expressions for the torsion and
curvature tensor fields:
\begin{align}
\label{eq:torsion:coordinates}
T^{ij}_k&=\Gamma^{ij}_k-\Gamma^{ji}_k-\frac{\partial \pi^{ij}}{\partial x^k},\\
\label{eq:curvature:coordinates}
R^{ijk}_l&=\Gamma^{ir}_l\Gamma^{jk}_r-\Gamma^{jr}_l\Gamma^{ik}_r
+\pi^{ir}\frac{\partial\Gamma^{jk}_l}{\partial
x^r}-\pi^{jr}\frac{\partial\Gamma^{ik}_l}{\partial x^r}
-\frac{\partial \pi^{ij}}{\partial x^r}\Gamma^{rk}_l.
\end{align}

\begin{rem}
\label{rem:curvature:bianchi:vectorbundle}
Expressions (\ref{curvature}) and (\ref{eq:curvature:derivative})
remain valid for any contravariant connection on a vector bundle
$E$ provided we replace $\gamma$ by a section of $E$. In this case
Bianchi's identity (\ref{eq:Bianchi:1}) can be expressed as
\[ \bigodot_{\al_1,\al_2,\al_3} D_{\al_1}(R(\al_2,\al_3))-
\bigodot_{\al_1,\al_2,\al_3} R([\al_1,\al_2],\al_3)=0.\]
\end{rem}

If it happens that the contravariant connection is related to some
covariant connection by:
\[ \#D_\al\beta =\nabla_{\#\al}\#\beta,\]
(e.~g., if $D$ is induced by a covariant connection and $\Pi$ is parallel,
so $D_\al=\nabla_{\#\al}$ and $D\Pi=0$) the torsion and
curvature tensor fields are transformed by the musical homomorphism
to the usual torsion and tensor fields of $\nabla$:
\[ T^{\nabla}(\#\al,\#\beta)=\#T^{D}(\al,\beta),\qquad
R^{\nabla}(\#\al,\#\beta)\#\gamma=\#R^{D}(\al,\beta)\gamma.\]

\subsection{Geodesics}

For contravariant connections parallel transport can only
be defined along curves lying in symplectic leaves of $M$. The same
restriction applies to geodesics:

\begin{defn}
Let $(\gamma(t),\al(t))$ be a cotangent curve on $M$.
We say that $(\gamma,\al)$ is a \textsc{geodesic} if:
\begin{equation}
(D_\al \al)_{\gamma(t)}=0.
\end{equation}
\end{defn}

In local coordinates, a curve
$(\gamma(t),\al(t))=(x^1(t),\dots,x^m(t),\al_1(t),\dots,\al_m(t))$ is
a geodesic iff it satisfies the following system of ode's
\begin{equation}
\label{eq:geodesic:coordinates}
\left\{\begin{array}{ll}
\frac{dx^i(t)}{dt}=\pi^{ji}(x^1(t),\dots,x^m(t))\al_j(t),& \\
&(i=1,\dots,m)\\
\frac{d\al_i(t)}{dt}=-\Gamma^{jk}_i((x^1(t),\dots,x^m(t))\al_j\al_k.&
\end{array}\right.
\end{equation}
From this we have:

\begin{prop}
Let $M$ be a Poisson manifold, with a contravariant connection
$\Gamma$, and $x_0\in M$. Given $\al_{x_0}\in T^*_{x_0} M$, there is a unique
maximal geodesic $t\mapsto (\gamma(t),\al(t))$, starting at ${x_0}\in M$, with
$\al(0)=\al_{x_0}$.
\end{prop}

\begin{proof}
Choose a systems of coordinates $(x^1,\dots,x^m)$ centered at $x_0$.
By standard uniqueness and existence results for ode's, system
(\ref{eq:geodesic:coordinates}) has a unique solution such that
$(x^1(0),\dots,x^m(0),\al_1(0),\dots,\al_m(0))
=(0,\dots,0,\al_{x_{0,1}},\dots,\al_{x_{0,m}})$.
\end{proof}

The geodesic given by this proposition is called the geodesic
through $x_0$ with cotangent vector $\al_{x_0}$. Note that if $S$ is the
symplectic leaf through $x_0$ and $v\in T_{x_0} S$ is a vector tangent to $S$,
there can be several geodesics with this tangent vector at $x_0$.
However, for an $\F$-connection geodesics are uniquely determined
by tangent vectors and coincide with the geodesics of the covariant
connection induced on $S$.

The following result is the analogue of a well known result in
affine geometry:

\begin{prop}
Let $\Gamma$ be a contravariant connection on $M$. There exists a
unique contravariant connection on $M$ with the same geodesics and
zero torsion.
\end{prop}

\begin{proof}
Choose local coordinates on $M$ so $D$ has symbols $\Gamma^{ij}_k$,
and consider the set of functions
\begin{equation}
{^*\Gamma}^{ij}_k=\frac{1}{2}\left(\Gamma^{ij}_k+\Gamma^{ji}_k+
\frac{\partial \pi^{ij}}{\partial x^k}\right)
\end{equation}
One checks that if $\Gamma^{ij}_k$ and ${\tilde{\Gamma}}^{lm}_n$
are related by the transformation law (\ref{eq:symbols:transform}),
then ${^*\Gamma}^{ij}_k$ and ${^*\tilde{\Gamma}}_n^{lm}$ are also
related by the same transformation law. It follows that we have a
well defined contravariant connection $D^*$ on $M$. From the local
coordinate expressions for the torsion
(\ref{eq:torsion:coordinates}) and the geodesics
(\ref{eq:geodesic:coordinates}), we see that $D^*$ has zero torsion
and the same geodesics as $D$.

For uniqueness, let $D$ and $D^*$ be two connections with the same
geodesics and torsion 0. We let
\begin{equation}
S(\al,\beta)=D_\al \beta - D^*_\al \beta,\qquad
\al,\beta\in\Omega^1(M).
\end{equation}
Then $S$ is $C^\infty$-linear, so it is a tensor. Since the
connections have 0 torsion, we have:
\begin{align}
S(\al,\beta)-S(\beta,\al)&=(D_\al \beta -D_\beta \al)- (D^*_\al
\beta-D^*_\beta \al)\\
&=[\al,\beta]-[\al,\beta]=0.\notag
\end{align}
so $S$ is a symmetric tensor. Now if $\al_p\in T^*_p M$, we can
choose the geodesic (for $D$ and $D^*$) with cotangent vector
$\al_p$ and associated 1-form $\al$ along $\gamma$. We have
\begin{equation}
S(\al_p,\al_p)=D_\al \al-D^*_\al \al=0,
\end{equation}
so $S=0$ and $D=D^*$.
\end{proof}

\comment{For a contravariant connection, one can also define the
\textsc{exponential mapping} at $x$, denoted $\Exp_x$. It is the
mapping $\al_x\overset{\Exp_x}{\longmapsto}\gamma(1)$,
where $(\gamma,\al)$ is the
unique geodesic in $M$ with cotangent vector $\al_x$. The
exponential is a mapping $\Exp_x:V\to M$ defined in some
neighborhood $V$ of the origin in $T^*_x M$. As usual, we say that
the connection is \emph{complete} if geodesics are defined for all
$t$. In this case, the exponential is a globally defined map
$\Exp_x:T^*_x M\to M$. Note that, unlike the covariant case, the
exponential map is not a local diffeomorphism (except, of course, in the
symplectic case).}

\subsection{Poisson Connections}

Linear contravariant connections for which the Poisson tensor is parallel
play an important role. Recall that a covariant connection for
which the Poisson tensor is parallel exists iff the Poisson
manifold has constant rank (see e.~g.~\cite{Vaisman:book:1}, thm.~2.20).
On the other hand, for contravariant connections a simple argument
involving a partition of unity shows that we have:

\begin{prop}
Every Poisson manifold has a linear contravariant connection with
contravariant derivative $D$ such that $D\Pi=0$.
\end{prop}

\begin{proof}
Let $U_a$ be a domain of a chart $(x^1,\dots,x^m)$. On $U_a$, the
contravariant connection $D^{(a)}$ with symbols
\[\Gamma^{ij}_k=\frac{\partial \pi^{ij}}{\partial x^k}\]
satisfies $D^{(a)}\Pi=0$. If we take an open cover of $M$ by such
chart domains and if $\sum_a\phi^{(a)}=1$ is partition of unity
subordinated to this cover, then $D=\sum_a \phi^{(a)} D^{(a)}$ is a
connection on $M$ for which $\Pi$ is parallel.
\end{proof}

We shall call a contravariant connection on $M$ such that the
Poisson tensor $\Pi$ is parallel a \textsc{Poisson connection}.
In the symplectic case, these coincide with the symplectic
connections.

If a Poisson connection has vanishing torsion then it is an $\F$-connection:
since $D\Pi=0$, we have $D\#=\#D$, and from $T=0$ we conclude that for
$\al,\beta\in\Omega^1(M)$
\begin{align*}
\#\al=0 \Longrightarrow \#D_\al\beta&=\#D_\beta\al+\#[\al,\beta]\\
                                    &=D_\beta\#\al+[\#\al,\#\beta]=0.
\end{align*}
Therefore, a torsionless Poisson connection is in fact a family of connections
along the leaves of $M$: for each symplectic leaf $i:S\hookrightarrow
M$ there exists a unique covariant symplectic connection $\nabla^S$ on
$S$, such that
\[ i^*D_\al\beta=\nabla^S_{\#i^*\al}i^*\beta, \quad \al,\beta\in\Omega^1(M).\]
As we pointed out above, a non-regular Poisson manifold does not
admit covariant connections for which the Poisson tensor is
parallel. Therefore, in general, it is not possible to glue
together the covariant connections $\nabla^S$ to get a connection
on $M$. As the following example shows, the family form by these
connections will develop singularities at points where the rank
drops.

\begin{ex}
Consider a 2-dimensional non-abelian Lie algebra $\gg$ and choose a
basis $\set{\omega_1,\omega_2}$ such that:
\[ [\omega_1,\omega_2]=\omega_1.\]
On $\gg^*$ we take the Lie-Poisson bracket which relative to the
coordinates $(x^1,x^2)$ defined by the dual basis satisfies
$\set{x^1,x^2}=x^1$. Now consider the contravariant connection on
$\gg^*$ defined by:
\[ D_{dx^1} dx^1=D_{dx^2}dx^2=D_{dx^2}dx^1=0,\qquad D_{dx^1}dx^2=dx^2.\]
One checks easily that $D$ has zero torsion and $D\pi=0$. On the
other hand there is no globally defined covariant connection
$\nabla$ on $\gg^*$ such that $D_\al=\nabla_{\#\al}$. In fact, if
such a connection existed, then denoting by $\Gamma_{ij}^k$ its
Christoffel symbols, we should have
\[ \Gamma^{ij}_k=\pi^{il}\Gamma^j_{lk},\]
where $\Gamma^{ij}_k$ are the symbols of $D$. Taking $i=k=1$,
$j=2$, this would give
\[ 1=x^1\Gamma^{2}_{21},\]
which is impossible. Note that formally we obtain the solution
$\Gamma^{2}_{21}=\frac{1}{x^1}$, so there exists a singular connection
with singular set $x^1=0$. This is precisely the set of points where
the rank drops from 2 to 0.
\end{ex}

\section{Poisson Holonomy}

\subsection{Holonomy of a Symplectic Leaf} For a regular foliation the
topological behaviour close to a given leaf is controlled by the
holonomy of the leaf. For a singular foliation, as is the case of
the symplectic foliation of a Poisson manifold, there is in general
no such notion of holonomy (see, however, \cite{Dazord:article:1}
where holonomy is defined for transversely stable leaves). It turns
out that in the case of a Poisson manifold it is still possible to
introduce a notion of holonomy which also reflects the Poisson
geometry of nearby leaves. In this theory of holonomy,
contravariant connections play a significant role.

Let $M$ be a Poisson manifold and let $i:S\hookrightarrow M$ be a
symplectic leaf of $M$. Denote by $\nu(S)=T_SM/TS$ the normal
bundle to $S$ and by $p:\nu(S)\to S$ the natural projection.
By the tubular neighborhood theorem, there exists a smooth immersion
$\tilde{i}:\nu(S)\to M$ satisfying the following properties:
\begin{enumerate}
\item[i)] $\tilde{i}|_Z=i$, where $Z$ is the zero section of $\nu(S)$;
\item[ii)] $\tilde{i}$ maps the fibers of $\nu(S)$
transversely to the symplectic foliation of $M$;
\end{enumerate}

Assume that we have fixed such an immersion.
Each fiber $F_x=p^{-1}(x)$ determines a splitting
$T_x \nu(M)=T_x S\oplus T_x F_x$, so we have
a decomposition:
\begin{equation}
\label{eq:decompose:S}
T^*_x \nu(M)=T_x^* S\oplus T^*_u F_x, \qquad \text{where } (T_x
F_x)^0\simeq T^*_x S,\quad (T_x S)^0\simeq T^*_x F_x.
\end{equation}
Note that $T_x S=\text{Im}\#_x=\#(T_x F_x)^0$. For each $u\in F_x$
we have an analogous splitting $T_u\nu(M)=\#(T_u F_x)^0\oplus T_u F_x$,
so there is also a decomposition:
\begin{equation}
\label{eq:decompose:normal:bundle}
T^*_u \nu(M)=(T_u F_x)^0\oplus T^*_u F_x, \qquad \text{where }
T^*_u F_x\simeq (\#(T_u F_x)^0)^0.
\end{equation}

Each such immersion induces a unique Poisson structure on the total
space $\nu(S)$ such that $\tilde{i}:\nu(S)\to M$ is a Poisson map.
Also, on each fiber $F_x=p^{-1}(x)$ there is an induced
\emph{transverse Poisson structure} $\Pi_x^\perp$: The
corresponding bundle map $\#^\perp:T^*F_x\to T F_x$ is defined as
the composed map
\[
T^*F_x\overset{q^*_x}{\to}T^*_{F_x}\nu(S)
\overset{\#}{\to}T_{F_x}\nu(S)\overset{q_x}{\to}T F_x,
\]
where $q_x:T_{F_x}\nu(S)\to T F_x$ is the bundle projection from
the restricted tangent bundle $T_{F_x}\nu(S)$ onto $T F_x$
associated with the decomposition
(\ref{eq:decompose:normal:bundle}).

Now let $\al\in T_x M$. We decompose $\al$ according to
(\ref{eq:decompose:S}):
\[\al=\al^\parallel+\al^{\perp},\text{where } \al^\parallel\in (T_x
F_x)^0\simeq T^*_xS,\quad \al^{\perp}\in(T_x S)^0\simeq T^*_xF_x.\]
Since $F_x$ is a linear space, there is a natural identification
$T^*_x F_x\simeq T^*_u F_x$, and we denote by
$\tilde{\al}^\perp_u\in T^*_uF_x\simeq(\#(T_uF_x)^0)^0$ the
element corresponding to $\al^{\perp}$.
On the other hand, the composition of the musical isomorphism $\#$
with the differential of the projection $p:\nu(S)\to S$ induces an
isomorphism between the annihilator $(T_u F_x)^0$ and $T_x S$, so
we also have an isomorphism $(T_u F_x)^0\simeq T^*_x S$.
If we denote by $\tilde{\al}^\parallel_u\in (T_u F_x)^0$ the
element corresponding to $\al^{\parallel}$ under this isomorphism,
we have $p_*\#\tilde{\al}^\parallel=\#\al$.

Given a covector $\al\in T^*_xM$ we shall define its
\emph{horizontal lift} to $\nu(S)$ by
\[ h(u,\al)=\#\tilde{\al}^\parallel_u+\#^\perp\tilde{\al}^\perp_u\in T_u\nu(S).\]
By construction, we have property (CI)$^*$ of a contravariant
connection
\[ p_*h(u,\al)=\#\al,\quad u\in p^{-1}(x),\]
so this horizontal lift defines a kind of generalized contravariant
connection in $\nu(S)$. Note that it depends both on the immersion
and on the Poisson tensor.

Let $(\gamma(t),\al(t))$, $t\in [0,1]$, be a cotangent curve in the
symplectic leaf $S$ starting at $x=\gamma(0)$. If $u\in\nu(S)|_x$
is a point in the fiber over $x$, there exists an $\eps>0$ and a
horizontal curve $\tilde{\gamma}(t)$ in $\nu(S)$, defined for
$t\in[0,\eps)$, which satisfies:
\[
\left\{
\begin{array}{l}
\frac{d}{dt}\tilde{\gamma}(t)=h(\tilde{\gamma}(t),\al(t)), \qquad
t\in[0,\eps),\\ \\
\tilde{\gamma}(0)=u.
\end{array}
\right.\]
Moreover, we can choose a neighborhood $U_\gamma$ of
$0\in\nu(S)|_x$, such that for each $u\in U_\gamma$ the lift
$\tilde{\gamma}(t)$ with initial point $u$ is defined for all
$t\in[0,1]$.

If $(\gamma(t),\al(t))$ is a cotangent loop based at $x\in S$ then
this lift gives, by passing from initial to end point, a
diffeomorphism $H_S(\gamma,\al)$ of $U_\gamma$ into another
neighborhood $V_\gamma$ of $0\in\nu(S)|_x$, with the property that
$0$ is mapped to $0$. One extends the definition of $H_S$ for
piecewise smooth cotangent loops in the obvious way.

Denote by $\Agerm(F_x)$ the group of germs at $0$ of Poisson automorphisms
of $F_x$ which map $0$ to $0$.

\begin{prop}
\label{prop:Poisson:holonomy}
Let $(\gamma,\al),(\gamma',\al')$ be cotangent loops based at $x\in
S$, then:
\begin{enumerate}
\item[i)] $H_S(\gamma,\al)$ is an element of $\Agerm(F_x)$;
\item[ii)] $H_S((\gamma,\al)\cdot(\gamma',\al'))=H_S(\gamma,\al)\circ
H_S(\gamma',\al')$, where the dot denotes concatenation of
cotangent loops.
\end{enumerate}
\end{prop}

\begin{proof}
Let $(\gamma(t),\al(t))$ be a cotangent curve in $S$. For each $t$, we
have a trivialization of $p:\nu(S)\to S$ in a neighborhood of
$\gamma(t)$ such that $p(x,y)=x$. If $\al(t)=\sum a(t)
dx|_{\gamma(t)}+b(t)dy|_{\gamma(t)}$ we consider the 1-form with
constant coefficients
$\al_t=\sum a(t)dx+b(t)dy$. The lift of its restriction to $S$ defines
the time-dependent vector field:
\[ X_t=\#\tilde{\al}^{\parallel}_t+\#^{\perp}\tilde{\al}^{\perp}_t,\quad \text{where }
\tilde{\al}^{\parallel}_t\in (TF_{\gamma(t)})^0,\
\tilde{\al}^{\perp}_t\in T^*F_{\gamma(t)}\simeq(\#(TF_{\gamma(t)})^0)^0.\]
For each $t$, the transverse component $\tilde{\al}^{\perp}_t$ is a closed
1-form in $F_{\gamma(t)}$.

The lifts $\tilde{\gamma}$ of $\gamma$ are the integral curves of
the vector field $X_t$. We claim that the flow
$\phi^t$ of this vector field preserves the transverse Poisson
structure $\Pi^{\perp}$
\begin{equation}
\label{eq:holonomy:flow}
(\phi^{-t})_*\Pi^{\perp}_{\phi^t(u)}=\Pi^{\perp}_u,
\end{equation}
so (i) follows. Part (ii) also follows since we have just shown that we
can take $H_S(\gamma,\al)$ as the time-1 map of some flow.

To prove (\ref{eq:holonomy:flow}) we observe that
\[\frac{d}{dt} (\phi^{-t})_*\Pi^{\perp}_{\phi^t(u)}=(\phi^{-t})_*
{\left[\frac{d}{dh} (\phi^{-h})_*\Pi^{\perp}_{\phi^{h}(\phi^t(u))}\right]}_{h=0},\]
and we use the following lemma:
\begin{lem}
\label{lemma:Lie:derivative}
If $\al_1,\al_2\in T^*_uF_x\simeq(\#(T_uF_x)^0)^0$ then
\[{\left[\frac{d}{dh} (\phi^{-h})_*\Pi^{\perp}_{\phi^{h}(u)}\right]}_{h=0}(\al_1,\al_2)=
(\Lie_{X_t}\Pi)_u(\al_1,\al_2)\]
\end{lem}

Now we have
\begin{align*}
\Lie_{X_t}\Pi(\al_1,\al_2)&=\Lie_{\#\tilde{\al}^{\parallel}_t}\Pi(\al_1,\al_2)+
\Lie_{\#^{\perp}\tilde{\al}^{\perp}_t}\Pi(\al_1,\al_2)\\
&=\Lie_{\#\tilde{\al}^{\parallel}_t}\Pi(\al_1,\al_2)+
\Lie_{\#^{\perp}\tilde{\al}^{\perp}_t}\Pi^{\perp}(\al_1,\al_2)
\end{align*}
The transverse component vanishes since $\tilde{\al}^{\perp}_t$ is
a closed form in the fiber, for each $t$. For the parallel
component we write $\tilde{\al}^{\parallel}_t=\sum_i a_i dx^i$, and
we compute
\[ \Lie_{\#\tilde{\al}^{\parallel}_t}\Pi=\sum_i
\left(a_i\Lie_{\#dx^i}\Pi+\#da_i\wedge\#dx^i\right).\]
But $dx^i\in (TF_x)^0$ and since $\al_1,\al_2\in
(\#(T_uF_x)^0)^0$ we conclude that
\[ \Lie_{\#\tilde{\al}^{\parallel}_t}\Pi(\al_1,\al_2)=
\sum_i a_i\Lie_{\#dx^i}\Pi(\al_1,\al_2)=0,\]
so the parallel component also vanishes.

It remains to prove lemma \ref{lemma:Lie:derivative}. We note
that for any $\al\in T^*_uF_x$ we have
$q^*_{\phi^h(u)}(\phi^{-h})^*\al-(\phi^{-h})^*q_u^*\al\in(T
F_{p(\phi^h(u))})^0$. Using this remark we find:
\begin{align*}
&{\left[\frac{d}{dh}
(\phi^{-h})_*\Pi^{\perp}_{\phi^{h}(u)}\right]}_{h=0}(\al_1,\al_2)=\\
&=\lim_{h\to 0}\frac{1}{h}\left[\Pi_{\phi^h(u)}(q^*_{\phi^h(u)}(\phi^{-h})^*\al_1,
q^*_{\phi^h(u)}(\phi^{-h})^*\al_2)-\Pi_u(q_u^*\al_1,q_u^*\al_2)\right]\\
&=
\Pi_{\phi^h(u)}\left[\Pi_{\phi^h(u)}((\phi^{-h})^*q_u^*\al_1,(\phi^{-h})^*q_u^*\al_2)-
\Pi_u(q_u^*\al_1,q_u^*\al_2)\right]\\
&=(\Lie_{X_t}\Pi)_u(q^*_u\al_1,q^*_u\al_2),
\end{align*}
so the lemma follows.
\end{proof}

Denoting by $\Omega_*(S,x)$ the group of piecewise smooth cotangent
loops, we see that we have a group homomorphism
$H_S:\Omega_*(S,x)\to\Agerm(F_x)$, which will be called the
\emph{Poisson holonomy homomorphism} of the leaf $S$. This Poisson
holonomy homomorphism depends on the immersion
$\tilde{i}:\nu(S)\to M$, but two different
immersions lead to conjugate homomorphisms.

\begin{ex} Let $S$ be a regular leaf of a Poisson manifold $M$. In
decomposition \ref{eq:decompose:normal:bundle} we can identify
$(T_u F_x)^0\simeq T^*_u S_u$ and $(T_u S_u)^0\simeq T^*_u F_u$, where
$S_u$ is the symplectic leaf through $u$. It follows that the horizontal
lift $h(u,\al)$ is the unique tangent vector in $T_u S_u$ which projects to
$\#\al$. We conclude that for a regular leaf the Poisson holonomy
coincides with the usual holonomy.
\end{ex}

\begin{ex}
\label{ex:holonomy:not:invariant}
Let $\gg$ be some finite dimensional Lie algebra and consider on
$M=\gg^*$ the canonical linear Poisson bracket. For the singular
leaf $S=\set{0}$ we have $\nu(S)\simeq\gg^*$ with $p(u)\equiv 0$
and the decomposition \ref{eq:decompose:normal:bundle} collapses.
Given a covector $\al\in T^*_0\gg^*=\gg$ we find
$h(u,\al)=\#_u\al=\ad^*\al\cdot u$. It follows that for a constant
cotangent loop $(0,\al)$ in $S$ we have
$H_S(0,\al)=\Ad^*(\exp(\al))$, which of course is a Poisson
automorphism of $F_0\simeq\gg^*$.
\end{ex}

\subsection{Reduced Poisson Holonomy}
As example \ref{ex:holonomy:not:invariant} shows, Poisson holonomy is not a
homotopy invariant. Following the construction given in
\cite{Ginzburg:article:1} for the linear case, we can give a notion
of \emph{reduced Poisson holonomy} which is homotopy invariant.

For a Poisson manifold $M$ let us denote by $\Aut(M)$ the group of
Poisson diffeomorphisms of $M$, and by $\Aut^0(M)$ its connected component
of the identity: given $\phi\in\Aut^0(M)$ there exists a smooth family
$\phi_t\in\Aut(M)$, $t\in[0,1]$, such that $\phi_0=\text{id}$,
$\phi_1=\phi$, and $\phi_t$ is generated by a time-dependent vector field:
\[ \frac{d\phi_t}{dt}=X_t\circ\phi_t.\]
The vector field $X_t$ is an infinitesimal Poisson automorphism:
\[ \Lie_{X_t}\Pi=0.\]
We shall say that $\phi$ is a \emph{inner Poisson automorphism} or a
\emph{hamiltonian automorphism} if
there exists a smooth family of hamiltonian functions $h_t:M\to\Rr$ such
that $X_t=X_{h_t}=\#dh_t$. The set $\Inn(M)\subset \Aut(M)$ of inner Poisson
automorphisms is a normal subgroup, and we define the group of
\emph{outer Poisson automorphisms} of $M$ to be the quotient
$\Out(M)=\Aut(M)/\Inn(M)$

Recall that for a symplectic leaf $S$ we denote by
$\Agerm(F_x)$ the group of germs at $0$ of Poisson automorphisms
of $F_x$ which map $0$ to $0$. We shall also denote by $\Ogerm(F_x)$ the
corresponding group of germs of outer Poisson automorphisms.

\begin{prop}
\label{prop:reduced:holonomy}
Let $S$ be a symplectic leaf of $M$, with Poisson holonomy homomorphism
$H_S:\Omega_*(S,x)\to\Agerm(F_x)$. If
$(\gamma_1,\al_1)$ and $(\gamma_2,\al_2)$ are cotangent loops with
$\gamma_1\sim\gamma_2$ homotopic
then $H_S(\gamma_1,\al_1)$ and $H_S(\gamma_2,\al_2)$ represent the same
equivalence class in $\Ogerm(F_x)$.
\end{prop}

\begin{proof}
Since any piecewise smooth path $\gamma\subset S$ can be made into
a cotangent path, by property (ii) in proposition
\ref{prop:Poisson:holonomy} it is enough to show that for every
$x\in S$ there exists a neighborhood $U$ of $x$ in $S$ such that if
$\gamma(t)\subset U$ is a piecewise smooth loop based at $x$ and
$\al(t)\in T^*M$ is a piecewise smooth family with
$\#\al=\dot{\gamma}$ then $H_S(\gamma,\al)\in\Inn(F_x)$.

To see this we use the same notation as in the proof of proposition
\ref{prop:Poisson:holonomy}. In a trivializing neighborhood $U$ of
$p:\nu(S)\to S$ containing $x$, we can decompose the vector field
$X_t$ as:
\[ X_t=\#\tilde{\al}^{\parallel}_t+\#^{\perp}\tilde{\al}^{\perp}_t,\quad \text{where }
\tilde{\al}^{\parallel}_t\in (TF_{\gamma(t)})^0,\
\tilde{\al}^{\perp}_t\in T^*F_{\gamma(t)}\simeq(\#(TF_{\gamma(t)})^0)^0.\]
For each $t$, the transverse component $\tilde{\al}^{\perp}_t$ can
be taken a closed 1-form in $F_{\gamma(t)}$. It is clear that the
parallel component $\#\tilde{\al}^{\parallel}_t$ has no effect on
the holonomy. Hence we can assume that $S=\set{x}$, $F_x=M$,
$\gamma$ is a constant path and $\tilde{\al}^{\perp}_t=\al(t)$, so
\[ X_t=\#\al(t)=\#dh_t,\]
for some function $h_t$ defined in a neighborhood of $x$.
Since $H_S(\gamma,\al)$ is the time-1 flow of this hamiltonian vector field
we conclude that $H_S(\gamma,\al)\in \Inn(F_x)$.
\end{proof}

Given a loop $\gamma$ in $S$ we shall denote by
$\bar{H}_S(\gamma)\in\Ogerm(F_x)$ the equivalence class of
$H_S(\gamma,\al)$ for some piece-wise smooth family $\al(t)$
with $\#\al(t)=\gamma(t)$. The map $\bar{H}_S:\Omega(S,x)\to\Ogerm(F_x)$
will be called the \emph{reduced Poisson holonomy homomorphism} of $S$.
This maps extends to continuous loops and, by a standard argument, it
induces a homomorphism $\bar{H}_S:\pi_1(S,x)\to\Ogerm(F_x)$ where
$\pi_1(S,x)$ is the fundamental group (the use of the same letter
to denote both these maps should not be the cause of any confusion).

\subsection{Stability}
The reduced Poisson holonomy of a leaf carries information on the
behaviour of the Poisson structure in a neighborhood of the leaf.
The simplest result in this direction can be obtained as follows:
let us call $S$
\emph{transversely stable} if the transverse Poisson manifold $N$ is
stable near $S\cap N$, i.~e., if $N$ has arbitrarily small
neighborhoods of $N\cap S$ which are invariant under all
hamiltonian automorphisms.

\begin{thm}
\label{thm:local:stability:I}
(Local Stability I) Let $S$ be a compact, transversely stable leaf,
with finite reduced holonomy. Then $S$ is stable, i.~e., $S$ has
arbitrarily small neighborhoods which are invariant under all
hamiltonian automorphisms. Moreover, each symplectic leaf of $M$
near $S$ is a bundle over $S$ whose fiber is a finite union of
symplectic leaves of the transverse Poisson structure.
\end{thm}

\begin{proof}
Assume first that $S$ has trivial reduced holonomy.
We fix an embedding $\tilde{i}:\nu(S)\to M$ as above and a base
point $x_0\in S$. Also, we choose a Riemannian metric on $S$.

By compactness of $S$, there exists a number $c>0$ such that every
point $x\in S$ can be connected to $x_0$ by a smooth cotangent path
of length $<c$. For some inner product on $\nu(S)|_{x_0}$, let
$D_\eps$ be the disk of radius $\eps$ centered at $0$. For each
$\eps>0$, there exists a neighborhood $U\subset D_\eps$ such that:
\begin{enumerate}
\item[i)] for any piecewise-smooth cotangent path in $S$,
starting at $x_0$, with length $\le 2c$ and for any $u\in U$,
there exists a lifting with initial point $u$;
\item[ii)] the lifting of any cotangent loop based at $x_0$
with initial point $u\in U$ has end point in $U$;
\item[iii)] $U$ is invariant under all hamiltonian automorphisms;
\end{enumerate}
In fact, let $(\gamma_1,\al_1),\dots,(\gamma_k,\al_k)$ be cotangent
loops such that $\gamma_1,\dots,\gamma_k$ are generators of
$\pi_1(S,x_0)$, and let $\phi_i$ be Poisson diffeomorphisms which
represent the germs $H_S(\gamma_i,\al_i)$. Since the reduced
holonomy is trivial, there is a neighborhood $U'$ of $0$ in
$F_{x_0}=\nu(S)|_{x_0}$ such that
$U\subset\text{domain}(\phi_1)\cap\cdots\cap\text{domain}(\phi_k)$,
and $\phi_i|U'\in\Inn(F_{x_0})$, for all i. Since $S$ is
transversely stable, we can choose a smaller neighborhood $U\subset
U'$ invariant under all hamiltonian automorphisms.

Given $x\in S$ and a cotangent path $(\gamma,\al)$ connecting $x_0$
to $x$, let us denote by $\sigma_{(\gamma,\al)}:U\to F_{x}$ the
diffeomorphism defined by lifting. It follows from i) and ii) above
that if $(\gamma',\al')$ is a cotangent path homotopic to
$(\gamma,\al)$ then
$\sigma_{(\gamma,\al)}(U)=\sigma_{(\gamma',\al')}(U)$. It follows
from iii) that $\sigma_{(\gamma,\al)}(U)$ is also invariant under
all hamiltonian automorphisms.

Let $V$ be a neighborhood of $S$ in $M$. There exists $\eps(x)>0$
such that for the corresponding $U_x\subset D_{\eps(x)}$ we have
$\sigma_{(\gamma,\al)}(U_x)\subset V\cap F_x$. By compactness of
$S$, we can choose $\eps>0$ (independent of $x\in S$) such that for
the corresponding $U\subset D_{\eps}$ we have
\[ \sigma_{(\gamma,\al)}(U)\subset V\cap F_x\]
Set
\[ V_0=\bigcup_{(\gamma,\al)}\sigma_{(\gamma,\al)}(U).\]
Then $V_0\subset V$ is a open neighborhood of $S$ which is
invariant under all hamiltonian automorphisms of $M$.

If $u,u'\in V_0$ are two points
in the same symplectic leaf such that $p(u)=p(u')=x$, then there is a path
$\tilde{\gamma}$ in this symplectic leaf connecting these two points. It
follows from the decomposition (\ref{eq:decompose:normal:bundle}) that
there exists a cotangent loop $(\gamma,\al)$ in $S$ such that
$\tilde{\gamma}$ is a horizontal lift of this loop. Thus $u'$ is the image
of $u$ by $H_S(\gamma,\al)$ which is a hamiltonian automorphism of
$V_0\cap F_x$. Therefore, $u$ and $u'$ lie in the same symplectic of
$V_0\cap F_x$. We conclude that each symplectic leaf of $M$ near $S$
is a bundle over $S$ whose fiber is a symplectic leaf of the transverse
Poisson structure.

Assume now that $S$ has finite reduced Poisson holonomy. We let
$q:\tilde{S}\to S$ be a finite covering space such that
$q_*\pi_1(\tilde{S})=\Ker\bar{H}_S\subset\pi_1(S)$. If we embed
$\nu(S)$ into $M$ as above, and let $\nu(\tilde{S})$ be the pull
back bundle of $\nu(S)$ over $\tilde{S}$, we have a unique Poisson
structure in $\nu(\tilde{S})$ such that the natural map
$\nu(\tilde{S})\to\nu(S)$ is a Poisson map. Moreover, the reduced
Poisson holonomy of $\nu(\tilde{S})$ along $\tilde{S}$ is trivial,
so we can apply the above argument to $\nu(\tilde{S})$ and the
theorem follows.
\end{proof}

\begin{rem}
If a leaf $S$ is transversely stable and $x\in S$, let $N$ denote a
stable neighborhood of $F_x$. For each cotangent path
$(\gamma,\al)$, the Poisson holonomy $H_S(\gamma,\al)$ induces a
homeomorphism of the orbit space of $N$, for the transverse Poisson
structure, mapping zero to zero. If $(\gamma_1,\al_1)$ and
$(\gamma_2,\al_2)$ are cotangent loops such that
$H_S(\gamma_1,\al_1)$ and $H_S(\gamma_2,\al_2)$ represent the same
class in $\Ogerm(F_x)$, then they induce the same germ of
homeomorphism of the orbit space mapping zero to zero. In
\cite{Dazord:article:1} holonomy of a general, transversely stable,
foliation is defined using germs of homeomorphisms of the orbit
space, which in the case of a Poisson manifold coincide with these
ones.
\end{rem}

\subsection{Strict Poisson Holonomy}
Another problem raised by the local splitting theorem
and related to stability is whether one has a global splitting of
an entire neighborhood of a leaf $S$. Note that if a neighborhood
$V$ of $S$ has a Poisson splitting $S\times N$ then projection to
the first factor is a Poisson map. This motivates the

\begin{defn}
\label{defn:Poisson:neighborhood}
Let $M$ be a Poisson manifold and $i:S\hookrightarrow M$ a
symplectic leaf of $M$. A \textsc{Poisson tubular neighborhood} of $S$
is a smooth immersion $\tilde{i}:\nu(S)\to M$ satisfying:
\begin{enumerate}
\item[i)] $\tilde{i}|_Z=i$, where $Z$ is the zero section of $\nu(S)$;
\item[ii)] $\tilde{i}$ maps the fibers of $\nu(S)$
transversely to the symplectic foliation of $M$;
\item[iii)] For the Poisson structure on $\nu(S)$ induced from
$\tilde{i}$, the canonical projection $p:\nu(S)\to S$ is a Poisson map;
\end{enumerate}
\end{defn}

Suppose $S$ admits a Poisson tubular neighborhood. Then the regular
distribution $\#(\Ker p_*)^0$ is integrable and $S$, identified with the
zero section, is an integral leaf of this distribution.
Hence, we can consider the holonomy of $S$ (in the usual sense)
as a leaf of the corresponding foliation. We call this the \emph{strict
Poisson holonomy} of $S$, and we denote by $\check{H}_S:\Omega(S,x)\to
\Dgerm(F_x)$ the associated holonomy map, where
$\Dgerm(F_x)$ denotes the group of germs of diffeomorphisms of $F_x$ which map
$0$ to $0$. Strict Poisson holonomy is related to
reduced Poisson holonomy as follows.

\begin{prop}
Assume $S$ admits a Poisson tubular neighborhood. The map
$\check{H}_S:\Omega(S,x)\to \Dgerm(F_x)$ has image inside $\Agerm(F_x)$
and the following diagram commutes:
\[
\xymatrix{
\Omega(S,x)\ar[r]^{\check{H}_S} \ar[dr]_{\bar{H}_S} & \Agerm(F_x)\ar[d]\\
 &\Ogerm(F_x)}
\]
\end{prop}

\begin{proof}
Fix a Poisson tubular neighborhood $p:\nu(S)\to S$ and consider the
gene\-ra\-li\-zed connection in $\nu(S)$ defined by the distribution
$\#(\Ker p_*)^0$. Given a loop $\gamma(t)$ in $S$ there exists a
family of closed forms $\al^S_t\in \Omega^1(S)$ such that
$\#\al^S_t(\gamma(t))=\dot{\gamma}(t)$. The horizontal lifts of
this loop are integral curves of the time-dependent vector field
\[ \check{X}_t=\#p^*\al^S_t.\]
Since $dp^*\al^S_t=p^*d\al^S_t=0$, this vector field is an
infinitesimal Poisson automorphism. We conclude that the holonomy
maps $\check{H}_S(\gamma)$ are Poisson automorphisms.

Moreover, in the notation of the proof of proposition
\ref{prop:reduced:holonomy}, we have
$\check{X_t}=\#\al_t^\parallel$. It follows that if $(\gamma,\al)$
is a cotangent loop in $M$ then $H_S(\gamma,\al)$ and
$\check{H}_S(\gamma)$ represent the same class in $\Out(F_x)$.
\end{proof}

We can now state and prove the following splitting result:

\begin{thm}
\label{thm:local:stability:II}
(Local Stability II)
Suppose $i:S\hookrightarrow M$ is a compact symplectic leaf of a
Poisson manifold $M$ which admits a Poisson tubular neighborhood.
Assume further that $S$ has finite strict Poisson holonomy
and let $q:\tilde{S}\to S$ be the finite covering corresponding to
$\Ker \check{H}_S\subset\pi_1(S,x)$. Then there is a 
neighborhood $V$ of $S$ and a finite covering Poisson map
$\phi:\tilde{S}\times N\to V$, where $N$ is a transverse Poisson manifold to
$S$. If $S$ is transversely stable, then we can chose $N$ and $V$ to be 
stable neighborhoods.
\end{thm}

\begin{proof}
By a standard homotopy lifting argument, as in the end of the proof
of theorem \ref{thm:local:stability:I}, it is enough to consider
the special case where the holonomy is trivial. We must then show
that there is a neighborhood $V$ of $S$ and a Poisson
diffeomorphism $\phi:S\times N\to V$, where $N$ is a transverse
Poisson manifold to $S$.

Again, we fix an embedding $\tilde{i}:\nu(S)\to M$ as above and a
base point $x_0\in S$. Also, we choose a Riemannian metric on $S$.
By compactness of $S$, there exists a number $c>0$ such that every
point $x\in S$ can be connected to $x_0$ by a smooth cotangent path
of length $<c$. For some inner product on $\nu(S)|_{x_0}$, let
$D_\eps$ be the disk of radius $\eps$ centered at $0$. There
exists an $\eps>0$ such that: for any piecewise-smooth cotangent
path in $S$, starting at $x_0$, with length $\le 2c$ and for any
$u\in D_\eps$, there exists a lifting with initial point $u$.
Moreover, by shrinking $\eps$ if necessary, we can assume that the
lifting of any cotangent loop based at $x_0$ with initial point $u$
also ends at $u$. In fact, let
$(\gamma_1,\al_1),\dots,(\gamma_k,\al_k)$ be cotangent loops such
that $\gamma_1,\dots,\gamma_k$ are generators of $\pi_1(S,x_0)$,
and let $\phi_i$ be Poisson diffeomorphisms which represent the
germs $\check{H}_S(\gamma_i,\al_i)$. Then, since the holonomy is
trivial by assumption, there is a neighborhood $U$ of $0$ in
$\nu(S)|_{x_0}$ such that
$U\subset\text{domain}(\phi_1)\cap\cdots\cap\text{domain}(\phi_k)$,
and $\phi_i|U=$identity, for all i. We need only to choose $\eps$
such that $D_\eps\subset U$.

For each $u\in D_\eps$ we define a map $\sigma_u:S\to M$ as
follows: let $x\in S$ and connect $x$ to $x_0$ by a cotangent path
$(\gamma,\al)$ of length $<c$. Let $\tilde{\gamma}$ be the unique
lift of $(\gamma,\al)$ starting at $u$, and define
$\sigma_u(x)=\tilde{\gamma}(1)$. This map is well defined because
the holonomy is trivial. Also, $\sigma_u$ is clearly a local
embedding since $p\circ\sigma_u=$identity on $S$. Since $S$ is
compact we conclude that $\sigma_u$ is an embedding.

The map $\sigma_u$ clearly depends smoothly on $u$, and since the
holonomy is trivial, the map $u\mapsto\sigma_u(x)$, for a fixed $x$,
is one-to-one. It follows that the map $\phi:S\times D_\eps\to M$ given
by $(x,u)\mapsto \sigma_u(x)$ is a diffeomorphism onto a neighborhood
$V$ of $S$.

By hypothesis, $p:\nu(S)\to S$ is a Poisson map. On the other hand,
the composition $p\circ\phi:S\times D_\eps\to S$ is just
projection into the first factor, which is also a Poisson map. Then
$\phi$ must also be a Poisson map.

Finally, if $S$ is transversely stable, we can choose an open set 
$N\subset D_\eps$
stable for the transverse structure, so $V=\phi(S\times N)$ is a stable
neighborhood.
\end{proof}

For simply connected leaves we obtain:

\begin{cor}
Let $i:S\hookrightarrow M$ be a compact, simply connected, symplectic leaf
of a Poisson manifold $M$, which admits a Poisson tubular neighborhood.
Then there is a neighborhood $V$ of $S$ and a Poisson
diffeomorphism $\phi:S\times N\to V$, where $N$ is a transverse
Poisson manifold to $S$. If $S$ is transversely stable, then we can chose $N$ 
and $V$ to be stable neigborhoods.
\end{cor}

One should note that, in general, a leaf does not have a Poisson
tubular neighborhood, and so strict Poisson holonomy is not
defined. In the following example we give a Poisson manifold $M$
with a compact, simply connected, symplectic leaf $S$, which has no
Poisson tubular neighborhood. In particular, $M$ does not split as
$S\times N$ in a neighborhood of $S$.

\begin{ex}
First observe that $CP(n)$ is a coadjoint orbit of $U(n+1)$,
since the standard action of $U(n+1)$ on $CP(n)$ is a transitive hamiltonian
action. In fact, a theorem of Kostant says that, for a
compact Lie group, all hamiltonian $G$-spaces on which $G$ acts
transitively are coadjoint orbits. The argument goes as follows:
Let $\Phi:CP(n)\to \mathfrak{u}^*(n+1)$ be the (equivariant) moment
map. Then $U(n+1)$ acts transitively on the image $Y=\Phi(CP(n))$, which
therefore is a coadjoint orbit. In fact, $\Phi:CP(n)\to Y$ is a
symplectomorphism and, since $CP(n)$ is compact, $\Phi$ is a covering map.
However, every coadjoint orbit of a compact Lie group is simply connected
(see \cite{Guillemin:book:1}, sect.~9.4), so this map is actually a
diffeomorphism.

Consider in particular the case $n=2$. We claim that $CP(2)$ is a
symplectic leaf of $\mathfrak{u}^*(3)$ which has no Poisson tubular
neighborhood. In fact, if $CP(2)$ had such a Poisson tubular
neighborhood then it would have trivial strict Poisson holonomy
and, by theorem \ref{thm:local:stability:II}, its normal bundle
$\nu(CP(n))$ would be trivial. But this is not the case, as can be
seen from the following standard argument: the total Chern class of
$CP(2)$ is $c=(1+a)^3=1+3a+3a^2$, where $a$ is a generator of
$H^2(CP(2),\Zz)$. The total Stiefel-Whitney class $w$ of $CP(2)$ is
the image of $c$ by the canonical homomorphism $H^2(CP(2),\Zz)\to
H^2(CP(2),\Zz_2)$ and hence is non-zero. The total Stiefel-Whitney
class of the normal bundle $\nu(CP(2))$ is $w^{-1}$, which is
non-trivial. We conclude that $\nu(CP(2))$ is non-trivial.
\end{ex}

\subsection{Linear Poisson Holonomy}
Let $M$ be a Poisson manifold and $i:S\hookrightarrow M$ a
symplectic leaf of $M$ with Poisson holonomy homomorphism
$H_S:\Omega_*(S,x)\to\Agerm(F_x)$ (once a tubular neighborhood as
been fixed).

On $T_0 F_x\simeq F_x$ we consider the Poisson
bivector field $\Pi^L$ which is the linear approximation at $0$ to
the Poisson bracket on $F_x$. Also, we denote by $\Aut(F_x)$ the
set of linear Poisson automorphisms of $(F_x,\Pi^L)$. There is a
map $d:\Agerm(F_x)\to \Aut(F_x)$ which assigns to a germ of a
Poisson diffeomorphism of $(F_x,\Pi^\perp)$, mapping zero to zero,
its linear approximation.

\begin{defn}
The \textsc{linear Poisson holonomy} of the leaf $S$ is the homomorphism
$H_S^L\equiv dH_S:\Omega_*(S,x)\to \Aut(F_x)$.
\end{defn}

One can check that this notion of linear Poisson holonomy is essentially
the same as the one introduced in \cite{Ginzburg:article:1}.

To define the \emph{reduced linear Poisson holonomy} of the leaf $S$ one
can either show that the class of $H_S^L(\gamma,\al)$ in
$\Out(F_x)=\Aut(F_x)/\Inn(F_x)$ is homotopy invariant, or else take
the composition
$\bar{H}_S^L\equiv\bar{d}\bar{H}_S:\pi_1(S,x)\to\Out(F_x)$,
where $\bar{d}:\Ogerm(F_x)\to \Out(F_x)$ is the natural map.
Similarly, if $S$ admits a Poisson tubular neighborhood, one can define
the \emph{strict linear Poisson holonomy} has the composition
$\check{H}_S^L\equiv d\check{H}_S:\pi_1(S,x)\to GL(F_x)$.

One can give a differential operator formulation for linear Poisson
holonomy similar to the Bott connection of ordinary foliation
theory. Instead of working with the normal bundle $\nu(S)=T_SM/TS$
it is convenient to use the dual bundle $\nu^*(S)$, also called the
conormal bundle. We have natural identifications
\[ \nu^*(S)=(\Ker\#)|_S=(T S)^0.\]
On $\nu^*(S)$ we have the following contravariant analogue of the Bott
connection:
Given a covector $\al\in T^*_SM$ and a section $\beta$ of $\nu^*(S)$,
take forms $\tilde{\al},\tilde{\beta}\in\Omega^1(M)$ such that
$\tilde{\al}_x=\al$, $\tilde{\beta}|_S=\beta$, and we set:
\begin{equation}
\label{Bott:connection}
D^S_\al\beta\equiv {[\tilde{\al},\tilde{\beta}]}|_S.
\end{equation}

To check that this definition is independent of the extensions considered,
we note that, by (\ref{Lie:bracket:forms}), it can also be written as
\begin{equation}
\label{Bott:connection:Lie}
D^S_\al\beta={\Lie_{\#\tilde{\al}}\tilde{\beta}}|_S.
\end{equation}
Expression (\ref{Bott:connection}) also shows that $D^S_\al \beta$ is in
the kernel of $\#$ and so is a section of $\nu^*(S)$.
Therefore, $D^S$ associates to each 1-form $\al$ on $M$ along $S$ a
differential operator $D_\al:\Gamma(\nu^*(S))\to\Gamma(\nu^*(S))$.

It is also easy to check that $D^S$ satisfies the analogue of properties
i)-iv) of proposition \ref{prop:derivative:properties:2}. Note however that,
in general, $D^S$ \emph{does not give} a contravariant connection in
$\nu^*(S)$, since it is defined only for 1-forms in $M$ along $S$. One
can now define parallel transport of fibers of $\nu^*(S)$ along cotangent
curves in $S$, and hence linear holonomy of $D^S$. The holonomy of $D^S$
coincides with the linear Poisson holonomy introduced above.

It is convenient to consider the connections $D^S$ all together, rather than
leaf by leaf, so we set:

\begin{defn}
A linear contravariant connection $D$ on $M$ is called a
\textsc{basic connection} if
\begin{enumerate}
\item[i)] $D$ restricts to $D^S$ on each leaf $S$, i.~e.,
if $\al,\beta\in\Omega^1(M)$ and $\#\beta|_S=0$ then
\[ D_\al\beta |_S=D^S_\al\beta.\]
\item[ii)] $D$ preserves the Poisson tensor, i.~e.,
\[D\Pi=0.\]
\end{enumerate}
\end{defn}

It is clear that one can also define linear Poisson holonomy starting
with some basic connection. The holonomy of this basic connection
determines maps of each cotangent space $T^*_x M$ which map $\ker\#_x$
isomorphically into itself, and these are the linear Poisson holonomy maps.

Basic connections always exist:

\begin{prop}
Every Poisson manifold has basic connections.
If $D$ is a basic connection with curvature tensor $R$, and $\gamma$ is
a 1-form such that $\#\gamma|_S=0$, then
\[ R(\al,\beta)\gamma|_S=0.\]
\end{prop}

\begin{proof}
Assume first that $M\simeq\Rr^m$, with coordinates $(x^1,\dots,x^m)$.
We define a contravariant connection on $M$ by setting
\[ D_{dx^j}\beta=[dx^j,\beta].\]
Then, obviously, if $S$ is a leaf of $M$ and $\#\beta|_S=0$ we have
\[ D_{dx^j}\beta |_S=D^S_{dx^j}\beta.\]
It follows that for any 1-form $\al$ we have
\[ D_\al\beta |_S=D^S_\al\beta.\]
Moreover, $D\Pi=0$ so $D$ is a basic connection.

For an arbitrary Poisson manifold $M$ we choose an open cover
$\set{U^{(a)}}$, with a partition of unity $\sum_a \phi_a=1$
subordinated to this cover, and such that on each $U^{(a)}$ there is a basic
connection $D^{(a)}$. Then $D=\sum_a \phi_a D^{(a)}$ is a
basic connection.

If $D$ is any basic connection and $\#\gamma|_S=0$, we have
$D_\al\gamma |_S=[\al,\gamma]|_S$ for any 1-form $\al$, so
expression (\ref{eq:curvature:derivative}) for the curvature tensor, gives
\[ R(\al,\beta)\gamma|_S=[\al,[\beta,\gamma]]|_S-[\beta,[\al,\gamma]]|_S-
[[\al,\beta],\gamma]|_S.\] But the right hand side is zero, because
of Jacobi identity.
\end{proof}

\begin{rem}
Although the curvature of a basic connection vanishes along
$\ker\#$, the holonomy along $\#$ need not be discrete (this is
because of the presence of an extra term in the holonomy theorem
\ref{thm:holonomy}). Hence, in general, linear Poisson holonomy is
not discrete and also not homotopy invariant (cf.~example
\ref{ex:holonomy:not:invariant}). However, if one can find a basic
connection which is an $\F$-connection, then Poisson holonomy is
discrete. Such is the case for a regular Poisson manifold, where
(linear) Poisson holonomy coincides with standard (linear)
holonomy.
\end{rem}

To finish this section we state the following result which by now should
be obvious.

\begin{prop}
Let $M$ be a Poisson manifold, and $S$ a symplectic leaf which
admits a transverse measure $\mu$ invariant under the hamiltonian
flow. Then, for every cotangent path $(\gamma,\al)$ in $S$
\[\det H_S^L(\gamma,\al)=1,\]
where the determinant is computed relative to $\mu$.
\end{prop}

This result also follows from a formula of Ginzburg and Golubev,
proved in \cite{Ginzburg:article:1}, which states that for
\emph{any} measure $\mu$ on $M$ one has
\begin{equation}
\label{eq:Ginzburg:Golubev}
 \det H_S^L(\gamma,\al)=\exp(\int_{(\gamma,\al)}v_{\mu}),
\end{equation}
where $v_{\mu}$ is the modular vector field of the measure $\mu$ (see
section \ref{section:moudlar:class}) and the determinant is computed relative
to the measure induced by $\mu$ on the transverse fiber.
This formula shows that there is a strong relationship between the modular
class and Poisson holonomy. In the next section we will introduce invariants
of a Poisson manifold which generalize the modular class,
and we will make this relationship more precise.

\section{Characteristic Classes}

\subsection{Poisson-Chern-Weil Homomorphism}

The usual Chern-Weil theory for characteristic classes extend to
contravariant connections, as was observed in \cite{Vaisman:art:1}.
We give here a short account since we shall need characteristic classes
later in the section.

Consider a principal G-bundle $p:P\to M$ over a Poisson manifold, and choose
some contravariant connection $\Gamma$ on $P$. Given any symmetric,
$\Ad(G)$-invariant, $k$-multilinear function
\[ P:\gg\times\cdots\times\gg\to\Rr\]
we can define a $2k$-vector field $\lambda(\Gamma)(P)$ on $M$ as follows. If
$U_j$ is a trivializing neighborhood, $x\in U_j$ and
$\al_1,\dots,\al_{2k}\in T^*_xM$ then we set
\begin{equation}
\label{eq:Chern-Weil:homomorphism}
\lambda(\Gamma)(P)(\al_1,\dots,\al_{2k})=\sum_{\sigma\in S_{2k}} (-1)^{\sigma}
P(\Xi_j(\al_{\sigma(1)},\al_{\sigma(2)}),\dots,\Xi_j(\al_{\sigma(2k-1)},\al_{\sigma(2k)})).
\end{equation}
By the transformation rule for the curvature bivector fields, this formula
actually defines a 2k-vector field on the whole of $M$.

\begin{prop}
For any symmetric, invariant, $k$-multilinear function $P$, the  $2k$-vector
field $\lambda(\Gamma)(P)$ is closed:
\begin{equation}
\delta \lambda(\Gamma)(P)=0.
\end{equation}
\end{prop}

\begin{proof}
We compute
\begin{align*}
\delta \lambda P(\Xi_j,\dots,\Xi_j)&=
k P(\delta \Xi_j,\dots,\Xi_j)\\
&=k P(\delta \Xi_j+[\Lambda_j,\Xi_j],\dots,\Xi_j)=0,
\end{align*}
where we have used first the linearity and symmetry of $P$, then the
$\Ad(G)$-invariance of $P$, and last the Bianchi identity.
\end{proof}
Therefore, to each invariant, symmetric, $k$-multilinear function
$P\in I^k(G)$ we can associate a  Poisson cohomology class
$[\lambda(\Gamma)(P)]\in H^{2k}_{\Pi}(M)$, and in fact we have:

\begin{prop}
\label{prop:connection:invariance:1}
The cohomology class $[\lambda(\Gamma)(P)]$ is independent of
the contravariant connection used to define it.
\end{prop}

\begin{proof}
Consider two contravariant connections $\Gamma^0$ and $\Gamma^1$ in $P$. Then
we have a family of connections $\Gamma^t$ with connection vector fields
$\Lambda^t_j=t\Lambda^1_j+(1-t)\Lambda^0_j$. We denote by $\Xi^t_j$
its curvature bivector fields.
Also, the difference $\Lambda^{1,0}_j=\Lambda^1_j-\Lambda^0_j$
is a $\gg$-valued vector field. By the transformation rule
(\ref{eq:transform:connection}), given $P\in I^k(G)$, we
get a well defined $(2k-1)$-vector field $\lambda(\Gamma^1,\Gamma^0)(P)$
by setting
\begin{multline}
\label{eq:invariants:2}
\lambda(\Gamma^1,\Gamma^0)(P)(\al_1,\dots,\al_{2k-1})=\\
k \sum_{\sigma\in S_{2k-1}} (-1)^{\sigma}\int_{0}^1
P(\Lambda^{1,0}_j(\al_{\sigma(1)}),\Xi^t_j(\al_{\sigma(2)},\al_{\sigma(3)}),
\dots,\Xi^t_j(\al_{\sigma(2k-2)},\al_{\sigma(2k-1)}))dt.
\end{multline}
We claim that
\begin{equation}
\label{eq:delta:lambda:2}
\delta\lambda(\Gamma^1,\Gamma^0)=\lambda(\Gamma^1)(P)-\lambda(\Gamma^0)(P),
\end{equation}
so $[\lambda(\Gamma^1)(P)]=[\lambda(\Gamma^0)(P)]$.

To prove (\ref{eq:delta:lambda:2}), we note that if we differentiate
the structure equation (\ref{eq:first:structure:equation}) we obtain
\begin{equation}
\label{eq:tderivative:curvature}
\frac{d}{dt}\Xi^t_j=\delta \Lambda^{1,0}_j+[\Lambda^t_j,\Lambda^{1,0}_j].
\end{equation}
Hence, using Bianchi's identity, we have
\begin{align*}
k\delta &\int_0^1 P(\Lambda^{1,0}_j,\Xi^t_j,\dots,\Xi^t_j)dt=\\
&=k\int_0^1 P(\delta\Lambda^{1,0}_j,\Xi^t_j,\dots,\Xi^t_j)
+P(\Lambda^{1,0}_j,\delta\Xi^t_j,\dots,\Xi^t_j)
+P(\Lambda^{1,0}_j,\Xi^t_j,\dots,\delta\Xi^t_j)dt\\
&=
k\int_0^1
P(\frac{d}{dt}\Xi^t_j-[\Lambda^t_j,\Lambda^{1,0}_j],\Xi^t_j,\dots,\Xi^t_j)\\
&\qquad \qquad \qquad \qquad
-P(\Lambda^{1,0}_j,[\Lambda^t_j,\Xi^t_j],\dots,\Xi^t_j)
-P(\Lambda^{1,0}_j,\Xi^t_j,\dots,[\Lambda^t_j,\Xi^t_j])dt\\
&=
k\int_0^1
P(\frac{d}{dt}\Xi^t_j,\Xi^t_j,\dots,\Xi^t_j)dt\\
&=
\int_0^1
\frac{d}{dt}P(\Xi^t_j,\Xi^t_j,\dots,\Xi^t_j)dt=
P(\Xi^1_j,\dots,\Xi^1_j)-P(\Xi^0_j,\dots,\Xi^0_j).
\end{align*}
\end{proof}

If we set
\[ I^*(G)=\bigoplus_{k\ge 0}I^k(G),\]
the assignment $P\mapsto [\lambda(\Gamma)(P)]$ gives a map
$I^*(G)\to H^*_\Pi(M)$. This map is in fact a ring
homomorphism.

\begin{prop}
\label{prop:Poisson:Chern:Weil}
The following diagram commutes
\[
\xymatrix{
I^*(G)\ar[r]\ar[dr]& H^*(M) \ar[d]^{\#} \\
& H^*_\Pi(M)}
\]
where on the top row we have the Chern-Weil homomorphism.
\end{prop}

\begin{proof}
Choose a contravariant connection $\Gamma$ in $P$ which is induced by a
covariant connection $\tilde{\Gamma}$. Given $P\in I^k(G)$, we
have a closed $(2k)$-form $\lambda(\tilde{\Gamma})(P)$ defined by a
formula analogous to (\ref{eq:Chern-Weil:homomorphism}), and which induces
the Chern-Weil homomorphism $I^*(G)\to H^*(M)$. We check easily that
\[ \#\lambda(\tilde{\Gamma})(P)=\lambda(\Gamma)(P),\]
so the proposition follows.
\end{proof}

Recall that the ring $I^*(GL_q(\Rr))$ is generated by elements
$P_k\in I^k(GL_q(\Rr))$ such that $P_k(A,\dots,A)=\sigma_k(A)$,
where $\set{\sigma_1,\dots,\sigma_q}$
are the elementary symmetric functions defined by:
\[ \det(\mu I-\frac{1}{2\pi}A)=\mu^q+\sigma_1(A)\mu^{q-1}+\cdots+\sigma_q(A).\]

Now consider a real vector bundle $p_E:E\to M$ over a Poisson manifold,
with fiber $F\simeq\Rr^q$ and let $p:P\to M$ be the associated principal
bundle with structure group $GL_q(\Rr)$. Choosing a contravariant connection
$\Gamma$ on $P$ one defines
the  \emph{kth Poisson-Pontrjagin class} of $E$ as
\[ p_k(E,\Pi)=[\lambda(\Gamma)(P_{2k})]\in H^{4k}_\Pi(M).\]
As usual, one does not need to consider the classes for odd $k$ since we have
\[ [\lambda(\Gamma)(P_{2k-1})]=0,\]
as can be seen by choosing a connection compatible with a riemannian metric.
It is clear from proposition \ref{prop:Poisson:Chern:Weil} that
\[ p_k(E,\Pi)=\#p_k(E).\]
where $p_k(E)$ are the standard Pontrjagin classes of $E$. Note also, that
if $r=\rank M=\max_{x\in M}(\rank\Pi_x)$ we have $p_k(E,\Pi)=0$ for $k>r/2$.

To compute these invariants one uses the contravariant derivative operator $D$
on $E$, associated with the contravariant connection $\Gamma$, and proceeds
as follows. For covectors $\al,\beta\in T_x M$, the curvature tensor $R$
defines a linear map $R_{\al,\beta}=R(\al,\beta):F_x\to F_x$ which satisfies
$R_{\al,\beta}=-R_{\beta,\al}$, and so $(\al,\beta)\to R_{\al,\beta}$ can be
considered as a $\gl(E)$-valued bivector field. By fixing a basis of local
sections, we have $F_x\simeq\Rr^q$ so we have $R_{\al,\beta}\in\gl_q(\Rr)$.
(this matrix representation of $R_{\al,\beta}$ is defined only up to a change
of basis in $\Rr^q$). Hence, if
\[ P:\gl_q(\Rr)\times\cdots\times\gl_q(\Rr)\to\Rr\]
is a symmetric, $k$-multilinear function, $\Ad(GL_q(\Rr))$-invariant, we
a have a $2k$-vector field $\lambda(R)(P)$ on $M$ defined by
\begin{equation}
\label{eq:Chern-Weil:homomorphism:connection}
\lambda(R)(P)(\al_1,\dots,\al_{2k})=\sum_{\sigma\in S_{2k}} (-1)^{\sigma}
P(R_{\al_{\sigma(1),\sigma(2)}},\dots,R_{\al_{\sigma(2k-1),\sigma(2k)}}).
\end{equation}
It is easy to see that $\lambda(\Gamma)(P)=\lambda(R)(P)$, so this gives
a procedure to compute the Poisson-Chern-Weil homomorphism and the
Poisson-Pontrjagin classes.

Similar considerations apply to other characteristic classes. One can
define, e.~g., the Poisson-Chern classes $c_k(E,\Pi)$ of a complex
vector bundle $E$ over a Poisson manifold, and they are just the images by
$\#$ of the usual Chern classes of $E$.

The fact that all these classes arise as image by $\#$ of some known classes
is perhaps a bit disappointing. However, we shall see below that
one can define Poisson secondary characteristic classes which are intrinsic of
Poisson geometry, and which do not arise as images by $\#$ of some de Rham
cohomology classes.

\subsection{Secondary Characteristic Classes}

We shall now introduce secondary characteristic classes of a
Poisson manifold. We will see that these classes give information on
the topology, as well as, the geometry of the symplectic foliation. As in the
theory of (regular) foliations, these classes appear when we compare two
connections, each from a distinguished class.

On the Poisson manifold $M$, with $\dim M=m$, we consider the
following data:
\begin{enumerate}
\item[i)] A basic connection $\Gamma^1$, with a contravariant derivative $D^1$;
\item[ii)] A linear contravariant connection $\Gamma^0$ induced by a riemannian
connection, so $D^0_\al=\nabla^0_{\#\al}$ with $\nabla g=0$ for some
riemannian metric $g$;
\end{enumerate}
Given an invariant, symmetric, $k$-multilinear function
$P\in I^k(GL(m,\Rr))$ we consider the $(2k-1)$-vector field
$\lambda(\Gamma^1,\Gamma^0)(P)$ given by (\ref{eq:invariants:2}).

\begin{prop}
\label{prop:invariants:2:closed}
If $k$ is odd, $\lambda(\Gamma^1,\Gamma^0)(P)$ is a closed $(2k-1)$-vector
field.
\end{prop}

\begin{proof}
According to (\ref{eq:delta:lambda:2}) we have
\[
\delta\lambda(\Gamma^1,\Gamma^0)=\lambda(\Gamma^1)(P)-\lambda(\Gamma^0)(P).
\]
and we claim that $\lambda(\Gamma^1)(P)=\lambda(\Gamma^0)(P)=0$ if $k$ is odd.

The proof that $\lambda(\Gamma^0)(P)=0$ is standard: since there exists a
metric such that $D^0 g=0$ we can reduce the structure group of $\Gamma^0$
to $O(m,\Rr)$, so the curvature bivector fields take there values in
$\mathfrak{so}(m,\Rr)$. But if $A\in \mathfrak{so}(m,\Rr)$, we have $P_k(A)$
for any elementary symmetric function, since $k$ is odd. Hence we obtain
$\lambda(\Gamma^0)(P)=0$.

Consider now the connection $\Gamma^1$. Given $x\in M$ we choose local
coordinates $(x^j,y^k)$ around $x$ as in the Weinstein splitting theorem:
\[ \Pi=\sum_{i=1}^{n}\frac{\partial}{\partial x^i}\wedge
\frac{\partial}{\partial x^{i+n}}+
\sum_{k,l}\phi_{kl}\frac{\partial}{\partial y^k}\wedge
\frac{\partial}{\partial y^l},\]
where $\phi_{kl}(x)=0$. Since $\Gamma^0$ is a basic connection,
we have:
\[ \Pi(D_\al dx^i,dx^j)=-\Pi(dx^i,D_\al dx^j),\qquad
R(\al,\beta)dy^k|_x=0.\]
It follows that $R(\al,\beta)_x$ is represented in the basis $(dx^j,dy^k)$
by a matrix of the form:
\begin{equation}
\label{eq:matrix:form}
\left(\begin{array}{cc}
B & 0 \\
C & 0
\end{array}\right),
\end{equation}
with $B$ a symplectic matrix. Now, if $A$ is any matrix of this form, it is
clear that $\det(\mu I-A)=\det(\mu I-\tilde{A})$, where $\tilde{A}$ is the
same as $A$ with $B=0$, i.~e., $\tilde{A}$ is symplectic. But if
$\tilde{A}$ is symplectic, we have $P_k(A)$
for any elementary symmetric function, since $k$ is odd.  Hence we obtain
also $\lambda(\Gamma^1)(P)_x=0$.
\end{proof}

Next we want to check that the Poisson cohomology class of
$\lambda(\Gamma^1,\Gamma^0)(P)$ is independent of the connections used to
define it.

Given connections $\Gamma^0,\Gamma^1,\Gamma^2$ we consider the family of
connections $\Gamma^{s,t}$ whose connection vector fields are
$\Lambda^{s,t}=(1-s-t)\Gamma^0+s\Gamma^1+t\Gamma^2$, where $(s,t)$ vary in
the standard 2-simplex $\Delta_2$. We introduce a $(2k-2)$-vector field
$\lambda(\Gamma^2,\Gamma^1,\Gamma^0)(P)$
given by a formula analogous to
(\ref{eq:Chern-Weil:homomorphism:connection}) and (\ref{eq:invariants:2}):
\begin{equation}
\label{eq:invariants:3}
\lambda(\Gamma^2,\Gamma^1,\Gamma^0)(P)=
k \sum_{\sigma\in S_{2k-2}} (-1)^{\sigma}\int_{\Delta_2}
P(\Lambda^{1,0}_j,\Lambda^{2,0}_j,\Xi^{s,t}_j,\dots,\Xi^{s,t}_j)dtds.
\end{equation}
Just like in the proof of proposition \ref{prop:connection:invariance:1}, one
shows that
\begin{equation}
\label{eq:delta:lambda:3}
\delta\lambda(\Gamma^2,\Gamma^1,\Gamma^0)=\lambda(\Gamma^1,\Gamma^0)(P)-
\lambda(\Gamma^2,\Gamma^0)(P)+\lambda(\Gamma^1,\Gamma^0)(P).
\end{equation}
Now, we can prove
\begin{prop}
\label{prop:independency:connection}
The Poisson cohomology class $[\lambda(\Gamma^1,\Gamma^0)(P)]$ is independent
of the connections used to define it.
\end{prop}

\begin{proof}
Let $\Gamma^1$ and $\tilde{\Gamma}^1$ (resp.~$\Gamma^0$ and
$\tilde{\Gamma}^0$) be basic connections (resp.~riemannian connections).
It follows from (\ref{eq:delta:lambda:3}) that
\begin{multline*}
\lambda(\Gamma^1,\Gamma^0)(P)-
\lambda(\tilde{\Gamma}^1,\tilde{\Gamma}^0)(P)=
\delta\lambda(\tilde{\Gamma}^1,\Gamma^0,\tilde{\Gamma}^0)(P)-
\delta\lambda(\Gamma^1,\tilde{\Gamma}^1,\Gamma^0)(P)\\
+\lambda(\tilde{\Gamma}^1,\Gamma^1)(P)
-\lambda(\Gamma^0,\tilde{\Gamma}^0)(P).
\end{multline*}
Hence, it is enough to show that the Poisson cohomology classes of
$\lambda(\tilde{\Gamma}^1,\Gamma^1)(P)$ and
$\lambda(\Gamma^0,\tilde{\Gamma}^0)(P)$ are trivial.

Consider first the basic connections $\tilde{\Gamma}^1$ and $\Gamma^1$.
The linear combination $\Gamma^{1,t}=(1-t)\Gamma^1+t\tilde{\Gamma}^1$ is also
a basic connection. If $x\in M$, we fix splitting coordinates $(x^j,y^k)$
around $x$ as in the proof of proposition \ref{prop:invariants:2:closed}.
Then we see that, with respect to the basis $\set{dx^j,dy^k}$, the matrix
representations of $D^1_\al$, $\tilde{D}^1_\al$ and $R^t(\al,\beta)$
are of the form (\ref{eq:matrix:form}). Hence, we conclude that if
$P\in I^k(GL(m,\Rr))$, with $k$ odd,
\[P(\tilde{D}^1_{\al_1}-D^1_{\al_1},R^t(\al_2,\al_3),\dots,
R^t(\al_{2k-2},\al_{2k-1}))=0.\]
Therefore, $\lambda(\tilde{\Gamma}^1,\Gamma^1)(P)=0$, whenever
$\tilde{\Gamma}^1$ and $\Gamma^1$ are basic connections.

Now consider the riemannian connections $\Gamma^0$ and
$\tilde{\Gamma}^0$. The linear combination
$\Gamma^{0,t}=(1-t)\tilde{\Gamma}^0+t\Gamma^0$ is also a riemannian
connection. All these connections are induced from covariant
riemannian connections $\nabla^0$, $\tilde{\nabla}^0$ and
$\nabla^{0,t}$, and we can define a $(2k-1)$-form
$\lambda(\nabla^0,\tilde{\nabla}^0)(P)$ by a formula analogous to
(\ref{eq:invariants:2}). Moreover, this form is closed (because $k$
is odd), and $
\#\lambda(\nabla^0,\tilde{\nabla}^0)(P)=\lambda(\Gamma^0,\tilde{\Gamma}^0)(P)$.
It follows from the homotopy invariance of $H^*(M)$, as in the
usual theory of characteristic classes of foliations (see
\cite{Bott:lectures:1}, page 29), that
\[ [\lambda(\nabla^0,\tilde{\nabla}^0)(P)]=[\lambda(\nabla^0,\nabla^0)(P)]=0.\]
Hence, the Poisson cohomology class $[\lambda(\nabla^0,\tilde{\nabla}^0)(P)]$
vanishes.

\end{proof}

\begin{rem}
The assumption that the riemannian connections are of the special form
$\nabla_{\#\al}$ was used in the proof to invoke the homotopy invariance of
$H^*(M)$. Poisson cohomology $H^*_\Pi(M)$ is not homotopy invariant, so
in defining the invariant $\lambda(\Gamma^1,\Gamma^0)(P)$ we cannot consider
an arbitrary
riemannian contravariant connection $\Gamma^0$. On the other hand, as we
pointed out above, in general a Poisson manifold does not admit a Poisson
connection of the form $\nabla_{\#\al}$. Hence, the basic connections are
``genuine'' contravariant connections, i.~e., not induced by any covariant
connection.
\end{rem}
\vskip 10pt

We define the \emph{secondary characteristic classes}
$\set{m_k(M)}$ of a Poisson manifold to be the Poisson cohomology
classes
\begin{equation}
m_k(M)=[\lambda(\Gamma^1,\Gamma^0)(P_{k})]\in H^{2k-1}_\Pi(M),\qquad
(k=1,3,\dots).
\end{equation}

If $M$ is a symplectic manifold then these classes obviously vanish. They
also vanish if $M=S\times N$ where $S$ is symplectic and $N$ has the zero
Poisson bracket. However, they
do not vanish for a general, regular, Poisson manifold
(see the examples below). Hence these
characteristic classes give information on both the Poisson geometry and the
topology of the symplectic foliation of $M$.
In the next section we give some explicit computations of these
classes, and in the following section we will show that the first
class coincides with the modular class of $M$ (up to a scalar factor).

\begin{rem}
In general, one can only define the characteristic classes $m_k$ for $k$ odd.
Assume, however, that $M$ admits flat riemannian connections and flat basic 
connections (we will see some non-trivial examples below). Then the proofs of 
propositions
\ref{prop:invariants:2:closed} and \ref{prop:independency:connection} can be
carried through, in the class of flat connections, for \emph{any} $k$.
Hence, in this case, one can define characteristic classes $m_k$ 
for \emph{any} $k$. 
\end{rem}

\subsection{Examples} We give a few types of Poisson manifolds where one
can compute some of the secondary characteristic classes.

\textbf{Euclidean spaces.}
Consider a Poisson manifold $M\simeq\Rr^m$, so we have global
coordinates $(x^1,\dots,x^m)$. To compute $\lambda(\Gamma^1,\Gamma^0)(P)$
we take as $\Gamma^0$ the flat connection determined by these global
coordinates, and as $\Gamma^1$ we take the basic connection defined by
\[ D_{dx^i} dx^j=[dx^i,dx^j]=\sum_{k}\frac{\partial\pi^{ij}}{\partial x^k}dx^k.\]
Since $P_1(A)=\frac{1}{2\pi}\tr(A)$, we find immediately that the
first characteristic class is
\begin{equation}
\label{eq:first:class}
m_1(M)=\frac{1}{2\pi}\sum_{i,j}
\frac{\partial\pi^{ij}}{\partial x^j}\frac{\partial}{\partial x^i}.
\end{equation}
To compute the second characteristic class, we note that
$D^t_{dx^i} dx^j=(1-t)D_{dx^i} dx^j$, and we compute its curvature:
\[ R^t(dx^i,dx^j)dx^k=-t(t-1) D_{[dx^i,dx^j]}dx^k.\]
Now, 
\begin{multline*} 
P_3(A,B,C)=\frac{1}{24\pi^3} \left[\tr(ABC)
-\frac{1}{2}(\tr A~\tr(BC)+\tr B~\tr(CA)+\tr C~\tr(AB))\right.\\
\left.-\frac{1}{2}\tr A~\tr B~\tr C \right]
\end{multline*}
and the expression for the characteristic class $m_3(M)$ is 
a certain homogeneous polynomial of degree $5$ involving the 
derivatives of order $\le 3$ of the components $\pi^{ij}$ of 
the Poisson tensor.
\vskip 15 pt

\textbf{Linear Poisson structures.}
Let $M=\gg^*$ with the  Lie-Poisson structure determined by the Lie algebra
$\gg$. Then, from the previous example, we see that the first
class is represented by the constant vector field
\[ m_1(\gg^*)(v)=\frac{1}{2\pi}\tr(\ad v).\]
In this case both the basic connection and the riemannian connection are
flat and so we can consider the classes $m_k$ for \emph{any} $k$. The 
computations simplify considerably, and see that all classes can be 
represented by constant multivector fields. For example, a straight forward
computation shows that
\begin{align*} 
m_2(\gg^*)(v_1,v_2,v_3)&=\frac{3!}{4\pi^2}K_2(v_1,[v_2,v_3]),\\
m_3(\gg^*)(v_1,\dots,v_5)&=\frac{1}{8\pi^3}
\sum_{\sigma\in S_5}
K_3(v_{\sigma(1)},[v_{\sigma(2)},v_{\sigma(3)}],[v_{\sigma(4)},v_{\sigma(5)}])
\end{align*}
where we have set
\[ K_j(v_1,\dots,v_j)=\tr(\ad v_1\cdots\ad v_j).\]
Note that $K_2$ is just the killing form. 

In this case it is possible to give a general formula for all characteristic
classes:
\begin{equation*}
m_k(\gg^*)(v_1,\dots,v_{2k-1})=\frac{1}{(2\pi)^k}\sum_{\sigma\in S_{2k-1}}
K_k(v_{\sigma(1)},[v_{\sigma(2)},v_{\sigma(3)}],\dots,
[v_{\sigma(2k-2)},v_{\sigma(2k-1)}])
\end{equation*}
The proof of these formulas involves a certain amount of computation using
Newton's identities for the elementary symmetric polynomials.

Incidentally, we note that the classes $m_k$ are $\ad$-invariant since each 
$K_j$ is an $\ad$-invariant multilinear form. Therefore, the classes 
$m_k(\gg^*)$ represent certain cohomology classes in the 
Lie algebra cohomology of $\gg$.
\vskip 15 pt

\textbf{Poisson-Lie Groups.}
Let $G$ be a connected Poisson-Lie group (see, e.~g., \cite{Lu:article:1}).
Then the set of left invariant 1-forms
$\Omega^1_{\text{Inv}}(G)$ is closed for the Lie bracket defined
by the Poisson bracket. Hence we can define a
basic connection $D^1$ in $G$ by requiring that
\begin{equation}
\label{eq:connection:group}
D_\al\beta=[\al,\beta], \qquad \forall \al,\beta\in\Omega^1_{\text{Inv}}(G).
\end{equation}
This connection is flat.

Let $D^0=\nabla_{\#\al}$ be the unique left invariant connection in
$G$ which for left invariant vector fields is given by
\[ \nabla_X Y=[X,Y],\qquad \forall X, Y\in\gg.\]
This connection is also flat.

We compute $\lambda(D^1,D^0)(P)$ and, generalizing the previous example, 
the classes $m_k(G)$ are all represented by the left invariant multivector 
fields:
\begin{equation*}
m_k(G)(\al_1,\dots,\al_{2k-1})=\frac{1}{(2\pi)^k}\sum_{\sigma\in S_{2k-1}}
K_k(\al_{\sigma(1)},[\al_{\sigma(2)},\al_{\sigma(3)}],\dots,
[\al_{\sigma(2k-2)},\al_{\sigma(2k-1)}])
\end{equation*}
where $\al_1,\dots,\al_n\in\Omega^1_{\text{Inv}}(G)$.
In these formulas, $[~,~]$, $\ad$ and $K_k$ are relative to the Lie algebra 
$\gg^*=\Omega^1_{\text{Inv}}(G)$.
\vskip 10 pt

\begin{rem}
Note that if the Poisson bracket in $G$ is not trivial, the contravariant 
connection
defined by (\ref{eq:connection:group}) is
\emph{not} left invariant, because left translation in the group is not a
Poisson map. These type of connections are studied in a complement to the
present paper, where we deal with invariant connections
(\cite{Fernandes:article:2}).
\end{rem}
\vskip 15 pt

\textbf{Regular Poisson manifolds.}
Let $M$ be a regular Poisson manifold of dimension $m$ and corank $q$.
First choose some riemannian connection determining a splitting
\[ T^*(M)=T^*(\mathcal{S})\oplus \nu^*(\mathcal{S}),\]
where $T^*(\mathcal{S})$ (resp.~$\nu^*(\mathcal{S})$) is the cotangent
(resp.~conormal) bundle to the symplectic foliation. We have a riemannian
connection $D^0$ such that:
\[ D^0_\al(\beta+\gamma)=\nabla_{\#\al}^{0,\parallel}\beta+
\nabla_{\#\al}^{0,\perp}\gamma, \]
where $\beta$ and $\gamma$, are sections of $T^*(\mathcal{S})$ and
$\nu^*(\mathcal{S})$, and $\nabla^{0,\parallel}$ and
$\nabla^{0,\perp}$,  are covariant riemannian connections in these bundles.

Becouse $M$ is regular, we can also choose a covariant connection
$\nabla^{1,\parallel}$ on $TM$
such that $\nabla^{1,\parallel}\Pi=0$. We define the basic connection $D^1$ on $M$
by setting
\[ D^1_\al(\beta+\gamma)=\nabla_{\#\al}^{1,\parallel}\beta+
\nabla_{\#\al}^{1,\perp}\gamma, \]
where $\nabla^{1,\perp}$ is a basic connection in $\nu(\mathcal{S})$ in
the usual sense of foliation theory (see \cite{Bott:lectures:1}, p. 33).
A computation shows that
\[
\lambda(D^1,D^0)(P)=\#\lambda(\nabla^{1,\perp},\nabla^{0,\perp})(\tilde{P}),
\]
where $\tilde{P}$ is the obvious restriction of $P\in
I^*(GL_m(\Rr))$ to $I^*(GL_q(\Rr))$.

It is well known in foliation theory (see \cite{Bott:lectures:1}, p.~66)
that the forms
\begin{align*}
c_k&=\lambda(\nabla^{1,\perp})(\tilde{P}_{k}),\qquad (1\le k\le q)\\
h_{2k-1}&=\lambda(\nabla^{1,\perp},\nabla^{0,\perp})(\tilde{P}_{2k-1}),
\qquad (1\le 2k-1\le q),
\end{align*}
satisfy
\begin{align}
\label{eq:Gelfand:Fuks}
dc_k&=0,\qquad (1\le k\le q)\\
dh_{2k-1}&=c_{2k-1},\qquad (1\le 2k-1\le q).
\end{align}
and so they can be used to define a homomorphism of graded algebras
\[ H^*(WO_q)\to H^*(M),\]
where $H^*(WO_q)$ is the relative Gelfand-Fuks cohomology of formal
vector fields in $\Rr^q$. This homomorphism is independent of the
connections and its image are the exotic or secondary characteristic
classes of foliation theory.

In this respect, the Poisson secondary characteristic classes are
simpler than the corresponding ones in foliation theory: the
$(2k-1)$-forms $\lambda(\nabla^{1,\perp},\nabla^{0,\perp})(\tilde{P}_{k})$
are not closed in general, but are closed along the symplectic leaves, so
its image under $\#$ is a closed $(2k-1)$-vector field and, hence,
define Poisson cohomology classes. Therefore, one has
\begin{equation}
m_{2k-1}(M)=[\#h_{2k-1}]
\end{equation}
but, in general, $m_{2k-1}$ is not in the image of $\#:H^*(M)\to H^*_{\Pi}(M)$.

Still, one can sometimes relate the two types of secondary
characteristic classes. Take, for example, the Godbillon-Vey class which
by definition is the cohomology class $w=[h_1c_1^q]\in H^{2q+1}(M)$
(it follows from relations (\ref{eq:Gelfand:Fuks}) that
$d(h_1c_1^q)=c_1^{q+1}=0$, so $h_1c_1^{q}$ does define a cohomology class).

\begin{prop}
If a regular Poisson manifold has a non-trivial
Godbillon-Vey class then it has a non-trivial first Poisson secondary
characteristic class.
\end{prop}

\begin{proof}
If $m_1(M)=[\#h_1]$ is trivial, we have $\#h_1=\#df$ for some smooth function
$f$, i.~ e., $h_1(\#\al)=df(\#\al)$. But $h_1$ is defined up to a 1-form in the
differential ideal that gives the symplectic foliation, so $h_1\wedge (dh_1)^ q=0$
and the the Godbillon-Vey class must vanish.
\end{proof}

On the other hand, it is perfectly possible for the Godbillon-Vey class to
vanish while $m_1(M)\not=0$. One such example is provided by the Reeb
foliation in $S^3$ with the leafwise area form (see \cite{Ginzburg:article:1,
Weinstein:article:2} for details on this example).

Another consequence of this relationship is that, for a regular Poisson
manifold $M$, the characteristic classes $m_k(M)=0$, for
$2k-1>q=\text{corank}(M)$.

As a special case, let us consider a Poisson manifold
of corank 1. The only non-vanishing class is $m_1(M)$. If the symplectic
foliation is transversely orientable, let $Z$ be a
trivializing section of the normal bundle. Let $\theta$ be the
corresponding 1-form that trivializes the conormal bundle. There
exists a 1-form $\eta$ such that
\[ d\theta=\eta\wedge\theta.\]
For $\nabla^{1,\perp}$ we choose a basic connection in
$\nu(\mathcal{S})$ such that
\[ \nabla^{1,\perp}_{X}Z=\eta(X)Z.\]
For $\nabla^{0,\perp}$ we choose a riemannian connection such
that
\[ \nabla^{0,\perp}_{X}Z=0.\]
These choices give
\[\lambda(\nabla^{1,\perp},\nabla^{0,\perp})(\tr)=\eta,\]
so we conclude that
\[ m_1(M)=\frac{1}{2\pi}[\#\eta].\]
In fact, in this case we have $h_1=\frac{1}{2\pi}\eta$ so
$w=\frac{1}{4\pi^2}\eta\wedge d\eta$ represents the Godbillon-Vey class.

If the symplectic foliation is not transversely orientable one can pass
to a double cover and apply the same reasoning.

\subsection{The Modular Class}
\label{section:moudlar:class}

The modular class of a Poisson manifold is an obstruction lying in
the first Poisson cohomology group $H^1_\Pi(M)$ to the existence of
a transverse invariant measure (see \cite{Weinstein:article:2} for
details on the modular class). It can be defined as follows: Let
$\mu$ be any measure in $M$ with associated divergence operator
$\text{div}_{\mu}X\equiv \Lie_X\mu/\mu$. Then one checks that the
map $f\mapsto \text{div}_\mu \#df$ is a derivation of $C^\infty(M)$
so defines a vector field $v_{\mu}$, called the
\emph{modular vector field} associated with the measure $\mu$. This
vector field is an infinitesimal automorphism of $\Pi$. If
$\mu'=a\mu$ is another measure, we have $v_{\mu'}=v_{\mu}+\#d\log a=
v_{\mu}+\delta\log a$, so in fact the modular class
\[ \mod(M)\equiv[v_\mu]\in H^1_{\Pi}(M)\]
is well defined and independent of $\mu$.

The examples in the previous section when compared to the computations
of the modular class done in \cite{Weinstein:article:2} suggest the
following

\begin{thm}
For any Poisson manifold $M$
\begin{equation}
\label{eq:invariant:modular:class}
m_1(M)=\frac{1}{2\pi}\text{mod~}(M).
\end{equation}
\end{thm}

\begin{proof}
Choose a basic connection $D^1$ and a riemannian connection $D^0$
relative to some metric on $M$. Let $\mu$ be the measure defined by
this metric. We claim that
\begin{equation}
\label{eq:invariant:modular:field}
\lambda(D^1,D^0)(\tr)=v_\mu,
\end{equation}
so (\ref{eq:invariant:modular:class}) follows.

Observe that it is enough to show that
(\ref{eq:invariant:modular:field}) holds on the regular points of
$M$, since the set of regular points is an open dense set and both
sides are smooth vector fields on $M$. So assume that $x\in M$ is a
regular point and pick Darboux coordinates $(x^1,\dots,x^m)$. If
$g=\left(\seq{dx^i,dx^j}\right)$ is the $m\times m$-matrix of inner
products of the $dx^i$'s, we have
\[ \mu=(\det g)^{\frac{1}{2}}dx^1\wedge\cdots\wedge dx^m.\]

As in the proofs of the previous section, relative to the basis
$\set{dx^1,\dots,dx^m}$, the operator $D^1_\al$ has a matrix
representation by a traceless matrix, so we only need to understand
what is the matrix representation, relative to this basis, of the riemannian
connection $D^0_\al=\nabla_{\#\al}$.

Since $\nabla$ is a metric connection, parallel transport preserves
the volume, and we have for any smooth function $f\in C^\infty(M)$:
\begin{align*}
0=&\nabla_{\#df}\mu\\
 =&\#df((\det g)^{\frac{1}{2}})dx^1\wedge\cdots\wedge dx^m+\\
 &\qquad \qquad \qquad +(\det g)^{\frac{1}{2}}(\nabla_{\#df}dx^1\wedge\cdots\wedge dx^m+
 \cdots+dx^1\wedge\cdots\wedge \nabla_{\#df}dx^m)\\
 =&\left(\#df((\det g)^{\frac{1}{2}})+(\det g)^{\frac{1}{2}}\tr \nabla_{\#df}\right)
 dx^1\wedge\cdots\wedge dx^m.
\end{align*}
So we conclude that:
\begin{equation}
\label{eq:modular:aux:1}
\tr (D^1_{df}-D^0_{df})\mu=\#df((\det g)^{\frac{1}{2}})dx^1\wedge\cdots\wedge dx^m.
\end{equation}

Now recall that $(x^1,\dots,x^m)$ were Darboux coordinates around a
regular point, so the form $dx^1\wedge\cdots\wedge dx^m$ is
preserved by the hamiltonian flows, and we have
\[ \Lie_{\#df}(dx^1\wedge\cdots\wedge dx^m)=0.\]
Hence, we conclude that:
\begin{equation}
\label{eq:modular:aux:2}
\Lie_{\#df}\mu=\#df((\det g)^{\frac{1}{2}})dx^1\wedge\cdots\wedge dx^m.
\end{equation}
Comparing (\ref{eq:modular:aux:1}) and (\ref{eq:modular:aux:2})
gives
\[ \tr (D^1_{df}-D^0_{df})=\text{div}_\mu \#df,\]
so relation (\ref{eq:invariant:modular:field}) holds.
\end{proof}

If $(\gamma(t),\al(t))$, $t\in[0,1]$, is a cotangent path and $X$ is a
vector field, one defines the integral
\[\int_{(\gamma,\al)}X=-\int_0^1 i_{X(\gamma(t))}\al(t) dt.\]
(For basic properties of this integral see
\cite{Ginzburg:article:1}). As a corollary of the theorem and the
Ginzburg and Golubev formula (\ref{eq:Ginzburg:Golubev}), we obtain:

\begin{cor}
Let $(\gamma,\al)$ be a cotangent loop in the symplectic leaf $S$.
Then
\begin{equation}
\det H_S^L(\gamma,\al)=\exp(\int_{(\gamma,\al)}tr(D^1-D^0)),
\end{equation}
where the determinant is relative to the transverse measure induced
by the volume element of the metric associated with $D^0$.
\end{cor}

\bibliographystyle{amsplain}
\def\lllll{}

\end{document}